\documentclass[11pt]{article}

\usepackage[letterpaper,margin=1in]{geometry}
\usepackage{amsmath,amssymb,amsthm}
\usepackage{graphicx}
\usepackage{booktabs} 
\usepackage{url}
\usepackage{array}
\usepackage{setspace}
\usepackage[authoryear,round]{natbib}
\usepackage{lineno}
\theoremstyle{definition}
\newtheorem{theorem}{Theorem}
\newtheorem{lemma}{Lemma}
\newtheorem{corollary}{Corollary}
\newtheorem{proposition}{Proposition}
\newtheorem{definition}{Definition}
\newtheorem{example}{Example}
\theoremstyle{remark}
\newtheorem{remark}{Remark}
\usepackage[svgnames]{xcolor}
\usepackage{soul}
\newcommand{\One}{\mathbf{1}}
\sethlcolor{PaleGreen}
\newif\ifmarked
\markedtrue   
\ifmarked
  
\else

\fi

 \usepackage{mdframed}
\ifmarked
  \newmdenv[hidealllines=true,leftline=true,linewidth=3pt,
            linecolor=yellow,innerleftmargin=6pt]{revblock}
\else
  
\fi

\newcommand{\R}{\mathbb{R}}
\newcommand{\norm}[1]{\left\|#1\right\|}

\title{Horseshoe-Type Shrinkage from Ratio-Generated Priors 
\thanks{The research of A. DasGupta was supported
in part by the National Science Foundation under
Grant DMS-2515247 (Purdue University).
The research of B. Vidakovic was supported
in part by the National Science Foundation under
Grant DMS-2515246 (Texas A\&M University).}
}

\author{
Radhika Kulkarni$^{1}$,
Anirban DasGupta$^{2}$,
and Brani Vidakovic$^{1}$%
\thanks{Corresponding author: Brani Vidakovic,
Department of Statistics, Texas A\&M University,
College Station, TX 77843, USA.
Email: \texttt{brani@tamu.edu}.}
\\[1em]
{\small $^{1}$Department of Statistics,
Texas A\&M University,
College Station, Texas, USA}\\
{\small $^{2}$Department of Statistics,
Purdue University,
West Lafayette, Indiana, USA}
}

\date{}

\begin{document}
\maketitle

\begin{abstract}
We develop a class of horseshoe-type shrinkage priors generated by
ratios of independent numerator and denominator random variables.
Large denominator values produce concentration near zero, small
denominator values generate heavy tails, and the numerator determines
scale and geometry and may also influence tail behavior. We establish
general origin and tail asymptotics, including power and logarithmic
singularities. 
Under a Gaussian likelihood and a spherically symmetric prior with
smoothly regularly varying tails, the posterior mean approaches the
observation for large signals, and its risk approaches
that of the unshrunk observation.
The denominator also induces a latent shrinkage factor. For a standard
Cauchy denominator, this factor has exactly the ${\rm Beta}(1/2,1/2)$
distribution, whatever the numerator. With a Gaussian numerator and
appropriate scale calibration, it recovers the classical horseshoe
shrinkage coefficient, whereas non-Gaussian numerators can change the
effective posterior shrinkage, even in its asymptotic rate.
The Cauchy/Cauchy prior illustrates the theory, with a logarithmic
singularity at the origin in one dimension, power singularities in
higher dimensions, and logarithmically modified Cauchy-order tails.
Its posterior mean shrinks smoothly, whereas its spatial median
behaves nearly as a threshold rule. We also compare the Cauchy/Cauchy posterior mean with fifteen established shrinkage procedures on the four standard Donoho--Johnstone test functions, finding competitive, though not uniformly superior, wavelet-denoising performance.
\end{abstract}


\section{Introduction}
\label{sec:introduction}

Sparse estimation is a central problem in modern statistics, signal
processing, and nonparametric function estimation. A canonical
formulation is the normal means model
\begin{eqnarray}
d_i
&=&
\theta_i+\sigma\varepsilon_i,
\qquad
\varepsilon_i\stackrel{\rm iid}{\sim}N(0,1),
\qquad
i=1,\ldots,n,
\label{eq:intro_normal_means}
\end{eqnarray}
where most of the unknown means $\theta_i$ are zero or close to zero,
and a comparatively small number carry appreciable signal. The same
model arises naturally in wavelet denoising. After an orthogonal
wavelet transform \citep{daubechies1992ten}, many smooth or spatially
inhomogeneous signals are represented by relatively few important
coefficients, while the remaining coefficients are dominated by noise.
This sparsity principle underlies classical wavelet thresholding and
shrinkage methods \citep{DonohoJohnstone1994,DonohoJohnstone1995}.

Bayesian approaches address the same problem through priors that
concentrate strongly near zero yet retain heavy tails, so that large
signals are not shrunk excessively. The horseshoe prior of
\citet{carvalho2010horseshoe} is an influential example. In its usual
global--local representation,
\begin{eqnarray}
\label{eq:intro_horseshoe}
\theta_i\mid\lambda_i,\tau
&\sim&
N(0,\tau^2\lambda_i^2),
\nonumber\\
\lambda_i
&\sim&
C^+(0,1),\\
\tau\mid\sigma &\sim& C^+(0, \sigma), \nonumber
\end{eqnarray}
where $\tau>0$ is a global scale, $\lambda_i$ is a coefficient-specific
local scale, $\sigma$ is the sampling standard deviation, and
$C^+(0,a)$ denotes the half-Cauchy distribution with density
$f(x)=\frac{2 a}{\pi (a^2 + x^2)} \One(x \geq 0)$.
Conditional on $\lambda_i$, the posterior mean under the Gaussian
likelihood has the familiar form
\begin{eqnarray}
E(\theta_i\mid d_i,\lambda_i,\tau,\sigma)
&=&
(1-\kappa_i)d_i,
\qquad
\kappa_i
=
{\sigma^2\over\sigma^2+\tau^2\lambda_i^2}.
\label{eq:intro_shrinkage_factor}
\end{eqnarray}
Small local scales produce strong shrinkage toward zero, whereas the
heavy tail of the half-Cauchy distribution allows large coefficients
to escape substantial attenuation. The horseshoe has become a standard
benchmark for continuous sparse shrinkage; see, among others,
\citet{polson2012local},
\citet{vanderpas2014horseshoe},
\citet{piironen2017sparsity}, and
\citet{bhadra2019lasso}.
More recent theoretical developments, methodological extensions, and
reviews include
\citet{vanderPas2021},
\citet{Bhadra2021},
\citet{KorobilisShimizu2022},
\citet{Banerjee2022},
\citet{KohnsSzendrei2024}, and
\citet{kollipara2026benchmarking}.
The broader literature on Bayesian shrinkage reaches well beyond the
horseshoe, including the three-parameter beta family
\citep{armagan2011generalized,armagan2013generalized} and, more
recently, priors induced by model fit and explained-variance
decompositions \citep{Zhang2022R2D2,AguilarBurkner2025,Yanchenko2026}.

Three features of the horseshoe are especially relevant here: its
marginal density is unbounded at the origin, its heavy tails protect
large signals from excessive shrinkage, and its latent shrinkage
coefficient has a U-shaped distribution, with substantial mass near
both zero and one. This paper asks where these features come from and
how far they extend.

We investigate them for priors generated by the ratio
\begin{eqnarray}
\Theta
&=&
\frac{X}{Y},
\label{eq:intro_ratio}
\end{eqnarray}
where $X\in\R^p$ and $Y\in\R$ are independent random variables with
densities $f$ and $g$, respectively. Quotient distributions have a long
history in probability theory \citep{Curtiss1941}. In the present
setting, large values of the denominator produce concentration at the
origin, and values near zero produce heavy tails. The numerator also
shapes the induced distribution, through its scale, its own
concentration and tails, and, when $p>1$, its directional geometry.

Our first contribution is a general theory of concentration at the
origin. We show that if $f$ is continuous and positive at zero, then
$E|Y|^p=\infty$ is sufficient for the density of $X/Y$ to diverge at
the origin, where $p$ is the dimension of $X$. Regular variation gives
sharper results: the right tail of $|Y|$ determines the strength of
the singularity, with explicit power behavior and a logarithmic
boundary that contains the scalar horseshoe as a special case.

Second, we establish how the distribution of the denominator near
zero, together with the tail of the numerator, determines the
polynomial tail order of the ratio, and how a borderline heavy
numerator can introduce a slowly varying correction without changing
that order. The analysis covers both scalar and multivariate settings.

Third, the ratio representation gives a natural multivariate
construction. A single scalar denominator acts simultaneously on all
coordinates of $X$. When the numerator has elliptical geometry, the
ratio preserves that geometry, while the denominator controls radial
concentration and tail thickness. This provides a shared-denominator
mechanism for multivariate and block shrinkage
\citep{CaiSilverman2001,AbramovichBesbeasSapatinas2002}.

Fourth, we study the inferential consequences of these heavy tails
under the Gaussian likelihood, including the large-signal behavior of
the posterior mean and its frequentist risk. These results link the
tails generated by the ratio construction directly to the robustness
required in sparse estimation.

Fifth, we investigate the shrinkage geometry induced by the
denominator. A Cauchy denominator generates the exact
$\operatorname{Beta}(1/2,1/2)$ latent geometry, independently of the
numerator. We show how this coordinate connects with classical
horseshoe shrinkage for Gaussian numerators, and how the relationship
changes, even in its asymptotic rate, for non-Gaussian numerators.

The theory is illustrated by representative ratio constructions,
several of which reproduce priors already in the literature, and in
detail by the Cauchy/Cauchy prior in both scalar and multivariate
settings. Under the Gaussian likelihood, we examine the corresponding
posterior mean and median shrinkage rules, using the spatial median in
the multivariate case.

The remainder of the paper is organized as follows.
Section~\ref{sec:horseshoe_review} reviews the classical horseshoe.
Section~\ref{sec:ratio_priors} develops ratio-generated priors and
establishes their origin singularities.
Section~\ref{sec:ratio_tails} investigates the tail behavior of
ratio-generated priors.
Section~\ref{sec:gaussian_sampling} studies posterior robustness under
the Gaussian likelihood.
Section~\ref{sec:kappa_geometry} examines denominator-induced
shrinkage geometry and its relationship to effective posterior
shrinkage.
Section~\ref{sec:cauchycauchy} presents representative constructions
and the Cauchy/Cauchy example.
Section~\ref{sec:dj_illustration} presents the comparative
wavelet-denoising illustration.
Section~\ref{sec:discussion} concludes.

\section{The Horseshoe Benchmark}
\label{sec:horseshoe_review}

The horseshoe hierarchy introduced in Section~\ref{sec:introduction}
provides a natural benchmark for the class developed below.  We focus on
the two features most relevant to the ratio construction: its shrinkage
behavior and the marginal shape of the induced prior.

With $\kappa_i$ defined in \eqref{eq:intro_shrinkage_factor},
conditioning on $\tau$ in \eqref{eq:intro_horseshoe} and taking the
standardized case $\sigma=\tau=1$ gives
$\kappa_i=(1+\lambda_i^2)^{-1}$, and the half-Cauchy distribution of
$\lambda_i$ induces
\begin{eqnarray*}
\kappa_i
&\sim&
{\rm Beta}\left({1\over2},{1\over2}\right).
\end{eqnarray*}
The U-shaped distribution of $\kappa_i$ places substantial probability mass near both
endpoints, producing the characteristic shrink-or-escape behavior of the
horseshoe \citep{carvalho2010horseshoe,polson2012local}.

\begin{figure}[h!tb]
\centering
\begin{minipage}{0.46\textwidth}
\centering
\includegraphics[width=\textwidth]{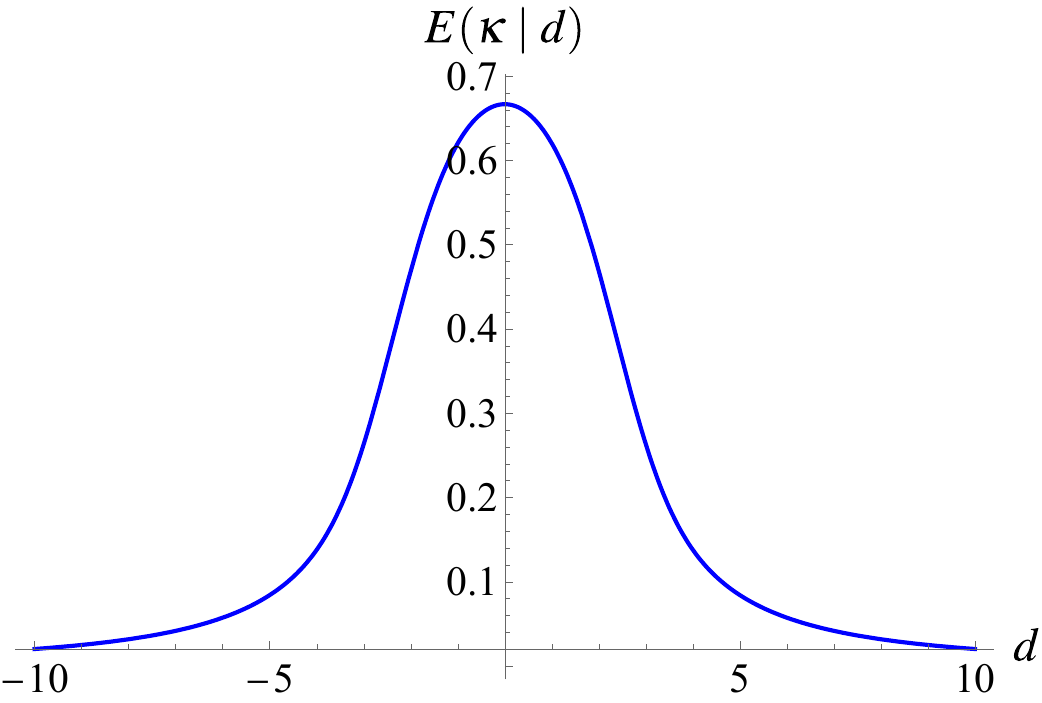}
\end{minipage}\hfill
\begin{minipage}{0.46\textwidth}
\centering
\includegraphics[width=\textwidth]{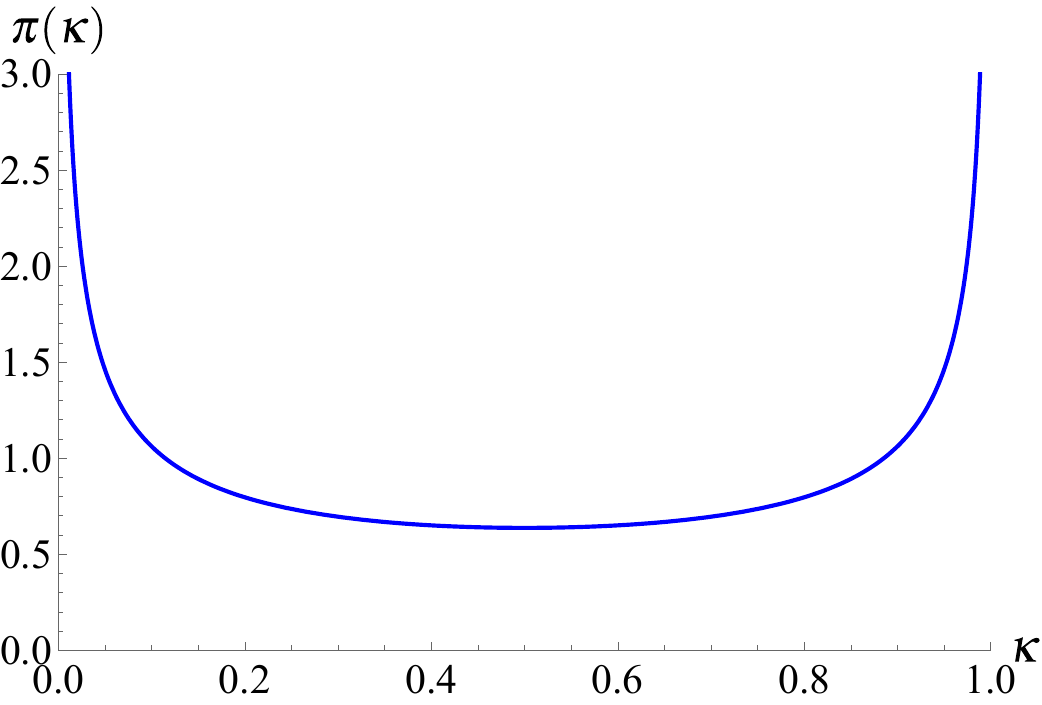}
\end{minipage}

\caption{Horseshoe shrinkage in the standardized normal means model.
Left: posterior expectation $E(\kappa\mid d)$ of the shrinkage coefficient.
Right: the horseshoe prior
$\kappa\sim{\rm Beta}(1/2,1/2)$, whose density is singular at both endpoints,
favoring either strong shrinkage, $\kappa$ near $1$, or essentially no
shrinkage, $\kappa$ near $0$.}

\label{fig:hs_kappa}
\end{figure}

In the standardized scalar model, the posterior density of $\kappa$ is
\begin{eqnarray*}
p(\kappa\mid d)
&\propto&
\exp\left(-{\kappa d^2\over2}\right)(1-\kappa)^{-1/2},
\qquad
0<\kappa<1,
\end{eqnarray*}
and
\begin{eqnarray*}
E(\theta\mid d)
&=&
\{1-E(\kappa\mid d)\}d.
\end{eqnarray*}
As $|d|$ increases, the posterior mass shifts toward smaller values of
$\kappa$, allowing large observations to escape substantial shrinkage.
Figure~\ref{fig:hs_kappa} displays these two aspects of the standardized
horseshoe construction.

The second benchmark feature is the marginal behavior of the horseshoe
prior for $\theta$, conditional on $\tau$. Its density exhibits a
logarithmic singularity at the origin and Cauchy-order tails, with the
following asymptotic equivalents:
\begin{align}
\pi_{\mathrm{HS}}(\theta\mid\tau)
&\sim
\frac{2}{\pi\sqrt{2\pi}\,\tau}\log\!\left(\frac{\tau}{|\theta|}\right),
~~~ \theta\to0,
\label{eq:hs_origin_order}\\
\pi_{\mathrm{HS}}(\theta\mid\tau)
&\sim
\frac{2\tau}{\pi\sqrt{2\pi}}\,|\theta|^{-2},
~~~ |\theta|\to\infty.
\label{eq:hs_tail_order}
\end{align}
 The
singularity provides strong concentration near zero, while the heavy tail
prevents persistent overshrinkage of large coefficients.  It is these
marginal properties, rather than the particular global--local hierarchy,
that motivate the ratio construction developed below.

An earlier related use of singular priors in wavelet shrinkage appears in
\citet{cutillo2008larger}, who considered a larger-posterior-mode rule
arising from the power-law prior
\begin{eqnarray*}
p(\lambda^2)
&\propto&
(\lambda^2)^{-k},
\qquad
\lambda^2>0.
\end{eqnarray*}
Under $\kappa=(1+\lambda^2)^{-1}$, the induced prior is
\begin{eqnarray*}
p(\kappa)
&\propto&
\kappa^{k-2}(1-\kappa)^{-k},
\qquad
0<\kappa<1.
\end{eqnarray*}
For $k=1$, this is proportional to
Haldane's improper prior
\citep{Haldane1932},
$\pi(\kappa) \propto \{\kappa(1-\kappa)\}^{-1}$ and is singular at both endpoints.  The
resulting posterior-mode rule yields wavelet thresholding and provides an
important precedent for exploiting singular prior behavior in sparse
estimation.

The horseshoe and the construction
in \citet{cutillo2008larger} demonstrate complementary manifestations of singular shrinkage
behavior.  The ratio construction developed below provides a direct
probabilistic mechanism linking concentration near zero, tail escape, and
the associated shrinkage geometry.

\section{Ratio-Generated Priors and Origin Singularities}
\label{sec:ratio_priors}

Let $X\in\R^p$ and $Y\in\R$ be independent random variables with
densities $f$ and $g$, respectively. Since both distributions are
absolutely continuous, the quotient $X/Y$ is well defined almost surely.

\begin{definition}[Ratio-generated prior]
\label{def:ratio_generated}
The distribution of
\begin{equation}
\Theta=\frac{X}{Y}
\label{eq:ratio_definition}
\end{equation}
is called a \emph{ratio-generated prior}.
\end{definition}

The definition deliberately imposes no symmetry, moment, or regularity
conditions; these enter only when particular properties of the induced
prior are studied. The change of variables $x=y\theta$ gives the joint
density $f_{\Theta,Y}(\theta,y)=|y|^pf(y\theta)g(y)$, and hence
\begin{equation}
\pi(\theta)=\int_{\R}|y|^pf(y\theta)g(y)\,dy,
\qquad \theta\in\R^p.
\label{eq:ratio_density}
\end{equation}
Every result in this section and the next is read off from this single
integral.

For shrinkage toward zero we impose the standing assumption that the
numerator is centrally symmetric, $f(x)=f(-x)$ for $x\in\R^p$. Then
\eqref{eq:ratio_density} immediately gives $\pi(-\theta)=\pi(\theta)$,
with no symmetry required of the denominator. Conditions on $f(0)$, on
the behavior of $g$ near zero and at infinity, and on moments of $X$
and $Y$ will be introduced only where they are needed.

We call a ratio-generated prior \emph{horseshoe-type} when it combines
an unbounded density at the origin with polynomially heavy tails. These
two properties are not built into the definition; the purpose of this
section and the next is to derive them from conditions on $X$ and $Y$.
Before doing so, we record a geometric consequence of
\eqref{eq:ratio_density} that requires no asymptotic analysis.

\begin{remark}
\label{rem:mvgeometry}
If the numerator has elliptical density
$$
f(x)=|\Sigma|^{-1/2}\psi(x^\top\Sigma^{-1}x),
$$
where $\Sigma$ is positive definite, then
$$
\pi(\theta)=|\Sigma|^{-1/2}\Psi(\theta^\top\Sigma^{-1}\theta),
\qquad
\Psi(s)=\int_{\R}|y|^p\psi(y^2s)g(y)\,dy.
$$
Division by a common scalar denominator therefore preserves the
elliptical geometry of the numerator while modifying its radial
concentration and tail behavior. The numerator controls the geometry,
and the denominator controls the spike--tail structure.
\end{remark}

The spike and the tail are governed by opposite ends of the
denominator. Under the standing symmetry assumption, let
$q(y)=g(y)+g(-y)$, $y>0$, denote the density of $|Y|$. For a unit vector
$e\in\R^p$ and $r>0$,
\begin{equation}
\pi(re)=\int_0^\infty y^pf(rye)q(y)\,dy.
\label{eq:radial_ratio_density}
\end{equation}
Small values of $r$ are reached through large values of $y$, so the
behavior of $q(y)$ as $y\to\infty$ determines the concentration of $\pi$
near the origin. Conversely, large values of $r$ require small $y$, so
the behavior of $q(y)$ as $y\downarrow0$ controls the far tail, subject
to the tail of the numerator. This section treats the origin;
Section~\ref{sec:ratio_tails} treats the tail.

We begin with a qualitative result. An unbounded density at the origin
requires only a moment condition on the denominator.

%
\begin{theorem}[Origin singularity]
\label{thm:origin_singularity}
Suppose that $f$ is continuous at the origin with $0<f(0)<\infty$. If
\begin{equation}
E|Y|^p=\infty,
\label{eq:infinite_p_moment}
\end{equation}
then $\pi(\theta)\to\infty$ as $\norm{\theta}\to0$.
\end{theorem}

\begin{proof}[Proof of Theorem~\ref{thm:origin_singularity}]
For each fixed $y$, $f(y\theta)\to f(0)$ as $\theta\to0$. Since the
integrand in \eqref{eq:ratio_density} is nonnegative, Fatou's lemma
gives
$$
\liminf_{\theta\to0}\pi(\theta)
\geq f(0)\int_{\R}|y|^pg(y)\,dy
=f(0)E|Y|^p=\infty.
$$
\end{proof}

\begin{remark}
Because $g$ is a density, any bounded set contributes finitely to
$E|Y|^p$, so \eqref{eq:infinite_p_moment} is purely a far-tail property
of the denominator. Theorem~\ref{thm:origin_singularity} thus formalizes
the contraction mechanism: occasional very large values of $|Y|$ pull
$X/Y$ toward zero often enough to make the density infinite there. For
a Gaussian numerator $X=\tau Z$ and $Y=1/\lambda$, the prior is the
normal scale mixture $\Theta\mid\lambda\sim N_p(0,\tau^2\lambda^2I_p)$,
and the condition becomes $E\lambda^{-p}=\infty$. For $p=1$ this
recovers the elementary fact that a normal scale mixture is unbounded at
zero exactly when $E\lambda^{-1}=\infty$; the theorem extends the
criterion to arbitrary numerators.
\end{remark}

Theorem~\ref{thm:origin_singularity} says whether a singularity occurs
but not how strong it is, and the strength matters: it determines how
much prior mass sits near zero, and hence how aggressively small
observations are shrunk. A quantitative answer requires a quantitative
description of the denominator's far tail. Suppose that
\begin{equation}
q(y)\sim c_Oy^{-1-\alpha}L_O(y),
\qquad y\to\infty,
\label{eq:q_origin_rv}
\end{equation}
where $c_O>0$, $\alpha>0$, and $L_O$ is slowly varying. We call
$\alpha$ the \emph{origin index} of the denominator. Since
$E|Y|^p<\infty$ when $\alpha>p$ and $E|Y|^p=\infty$ when $\alpha<p$,
the index $\alpha$ is compared with the dimension $p$, and three regimes
arise: $\alpha<p$, $\alpha=p$, and $\alpha>p$. Regular variation and
the associated integral asymptotics are developed in
\citet[Chapter VIII, Sections 8--9, pp.~275--283]{Feller1971}; for the
Potter bounds used in the proofs, see
\citet[Theorem~1.5.6, p.~25]{BinghamGoldieTeugels1987}.

In the first regime the denominator tail is heavy enough to produce a
power singularity.

%
\begin{theorem}[Power singularity at the origin]
\label{thm:origin_power}
Suppose that \eqref{eq:q_origin_rv} holds with $0<\alpha<p$. For a unit
vector $e\in\R^p$, assume that $u\mapsto f(ue)$ is bounded on $[0,1]$
and that, for some $\eta\in(0,p-\alpha)$,
$$
\int_1^\infty u^{p-\alpha-1+\eta}f(ue)\,du<\infty.
$$
Define
$$
C_\alpha(e)=\int_0^\infty u^{p-\alpha-1}f(ue)\,du.
$$
If $C_\alpha(e)>0$, then
\begin{equation}
\pi(re)\sim c_OC_\alpha(e)\,r^{\alpha-p}L_O(1/r),
\qquad r\downarrow0.
\label{eq:origin_power}
\end{equation}
\end{theorem}

The idea of the proof is a change of scale. Substituting $u=ry$ in
\eqref{eq:radial_ratio_density} turns the far tail of $q$ into the
factor $r^{\alpha-p}L_O(1/r)$, and Potter's bounds justify passing to
the limit under the integral; the details are in
Appendix~\ref{sec:prooftheorem2}.

The theorem separates the two roles cleanly. The exponent $\alpha-p$
and the slowly varying factor come from the denominator alone, so the
singularity is generated by the denominator even when $f(0)=0$. The
numerator enters only through the constant $C_\alpha(e)$, a weighted
integral of $f$ along the ray $\{ue:u>0\}$, which may depend on the
direction when $f$ is not spherically symmetric.

A divergent density need not mean that much probability sits near the
origin, so it is useful to translate \eqref{eq:origin_power} into prior
mass. Since $\norm{\Theta}\leq\varepsilon$ if and only if
$|Y|\geq\norm{X}/\varepsilon$, independence gives
$$
P(\norm{\Theta}\leq\varepsilon)
=E\left\{\bar Q\left(\frac{\norm{X}}{\varepsilon}\right)\right\},
\qquad \bar Q(y)=P(|Y|>y),
$$
the counterpart, at the origin, of the tail identity
\eqref{eq:tail_identity} of Section~\ref{sec:ratio_tails}. Under
\eqref{eq:q_origin_rv}, Karamata's theorem gives
$\bar Q(y)\sim(c_O/\alpha)y^{-\alpha}L_O(y)$. If
$E\norm{X}^{-\alpha-\eta}<\infty$ for some $\eta>0$, which holds when
$f$ is bounded near the origin and $\alpha<p$, the dominated-convergence
argument used for Theorem~\ref{thm:radial_tail}, with Potter bounds for
$\bar Q$ at infinity, gives
$$
P(\norm{\Theta}\leq\varepsilon)
\sim\frac{c_O}{\alpha}\,E\norm{X}^{-\alpha}\,
\varepsilon^\alpha L_O(1/\varepsilon),
\qquad \varepsilon\downarrow0.
$$
In polar coordinates,
$E\norm{X}^{-\alpha}
=\int_{S^{p-1}}C_\alpha(e)\,d\omega(e)$.
Thus the ball probability asymptotic is consistent
with the directional density asymptotic
\eqref{eq:origin_power}: the density diverges
like $r^{\alpha-p}$, whereas the probability
of a ball of radius $\varepsilon$ decays like
$\varepsilon^\alpha$. The smaller the
origin index, the more prior probability concentrates near zero.

As $\alpha$ increases to $p$, the exponent in \eqref{eq:origin_power}
vanishes and Theorem~\ref{thm:origin_power} no longer applies. This
boundary is the most important case in practice: for $p=1$, a Cauchy
denominator has $\alpha=1=p$, and so does the half-Cauchy local scale of
the horseshoe. At the boundary the singularity is slowly varying, and
its form is determined by the integrated slowly varying function below.

\begin{proposition}[Boundary case $\alpha=p$]
\label{prop:origin_log}
Suppose that \eqref{eq:q_origin_rv} holds with $\alpha=p$, where $L_O$
is positive and locally integrable, and define
$$
\widetilde L_O(x)=\int_1^x\frac{L_O(t)}{t}\,dt.
$$
Assume that $\widetilde L_O(x)\to\infty$, that $f$ is continuous at
the origin with $0<f(0)<\infty$, and that, for a unit vector
$e\in\R^p$, $u\mapsto f(ue)$ is locally bounded on $(0,\infty)$ and,
for some $\eta>0$,
\begin{equation}
\int_1^\infty u^{\eta-1}f(ue)\,du<\infty.
\label{eq:boundary_ray_integrability}
\end{equation}
Then
\begin{equation}
\pi(re)\sim c_Of(0)\widetilde L_O(1/r),
\qquad r\downarrow0.
\label{eq:origin_log_general}
\end{equation}
In particular, if $L_O(y)\to L_{O,0}\in(0,\infty)$, then
\begin{equation}
\pi(re)\sim c_Of(0)L_{O,0}\log(1/r),
\qquad r\downarrow0.
\label{eq:origin_log}
\end{equation}
\end{proposition}

Here the leading contribution comes from an intermediate range of
denominator values, $y_0<y<\delta/r$, over which $f(rye)$ is
essentially $f(0)$ and $q(y)$ contributes $L_O(y)/y$; integrating over
this range produces $\widetilde L_O(1/r)$. The two outer ranges are of
smaller order. The details are in Appendix~\ref{sec:proofproposition}.

The boundary differs from the power regime in one further respect: the
leading term in \eqref{eq:origin_log_general} does not depend on the
direction $e$. Only $f(0)$ enters, so the singularity is isotropic even
when the numerator is not spherically symmetric.

\begin{remark}[Summary of the origin behavior]
Assume that $0<f(0)<\infty$, together with the regularity conditions
of the respective results. The origin index $\alpha$, compared with the
dimension $p$, determines the behavior at zero: (i) if $0<\alpha<p$,
there is a power singularity of order $r^{\alpha-p}$
(Theorem~\ref{thm:origin_power}); (ii) if $\alpha=p$ and
$\widetilde L_O(x)\to\infty$, the singularity is slowly varying, and it
is logarithmic when $L_O$ approaches a positive constant
(Proposition~\ref{prop:origin_log}); (iii) if $\alpha>p$, or $\alpha=p$
with $\widetilde L_O(x)$ convergent, then $E|Y|^p<\infty$, and when, in
addition, $f$ is bounded, dominated convergence in
\eqref{eq:ratio_density} gives $\pi(\theta)\to f(0)E|Y|^p<\infty$ as
$\theta\to0$.
\end{remark}

The following example shows all of this in the case that motivates the
paper.

\begin{example}[Cauchy denominator]
Let $Y\sim C(0,1)$. Its folded density satisfies
$$
q(y)=\frac{2}{\pi(1+y^2)}\sim\frac{2}{\pi}y^{-2},
\qquad y\to\infty,
$$
so $\alpha=1$ and $L_O\equiv1$. For $p=1$,
Proposition~\ref{prop:origin_log} gives
$$
\pi(re)\sim\frac{2f(0)}{\pi}\log(1/r),
\qquad r\downarrow0.
$$
With $X\sim N(0,\tau^2)$, so that $f(0)=1/(\sqrt{2\pi}\,\tau)$, this is
exactly the horseshoe asymptotic \eqref{eq:hs_origin_order}: the folded
density of a half-Cauchy denominator is the same $q$, and $1/Y$ is the
half-Cauchy local scale. For $p>1$, the same denominator falls in the
power regime, and Theorem~\ref{thm:origin_power} gives
$$
\pi(re)\sim\frac{2C_1(e)}{\pi}r^{1-p},
\qquad r\downarrow0.
$$
The logarithmic spike of the scalar horseshoe is thus a boundary
phenomenon; with the same denominator, the spike becomes a power
singularity in every higher dimension.
\end{example}

\section{Tail Behavior of Ratio-Generated Priors}
\label{sec:ratio_tails}

Section~\ref{sec:ratio_priors} showed that the far tail of the
denominator governs concentration at the origin. We now turn to the
opposite end of $|Y|$: small denominators produce large values of
$\Theta$, and their frequency determines the tail of the prior. Let
$H(t)=P(|Y|\leq t)$, $t\geq0$. Since $\norm{\Theta}>r$ if and only if
$|Y|<\norm{X}/r$, independence gives the exact identity
\begin{equation}
P(\norm{\Theta}>r)=E\left\{H\left(\frac{\norm{X}}{r}\right)\right\}.
\label{eq:tail_identity}
\end{equation}

The radial tail of $\Theta$ is therefore an average of the distribution
function of $|Y|$ near zero. Equivalently, since
$\norm{\Theta}=\norm{X}\cdot|Y|^{-1}$ with independent factors, it is
the tail of a product, with $P(|Y|^{-1}>r)=H(1/r)$.

\begin{theorem}[Radial tails induced by small denominators]
\label{thm:radial_tail}
Suppose that the distribution function $H(t)=P(|Y|\leq t)$ satisfies
\begin{equation}
H(t)\sim c_Yt^\beta L_T(1/t),
\qquad t\downarrow0,
\label{eq:H_regular_variation}
\end{equation}
where $c_Y>0$, $\beta>0$, and $L_T$ is positive and slowly varying at
infinity. If $E\norm{X}^{\beta+\epsilon}<\infty$ for some $\epsilon>0$,
then
\begin{equation}
P(\norm{\Theta}>r)\sim c_YE\norm{X}^\beta\,r^{-\beta}L_T(r),
\qquad r\to\infty.
\label{eq:radial_tail}
\end{equation}
\end{theorem}

Theorem~\ref{thm:radial_tail} is Breiman's lemma applied to the product
$\norm{X}\cdot|Y|^{-1}$: by \eqref{eq:H_regular_variation},
$|Y|^{-1}$ is regularly varying with index $-\beta$, and the moment
condition on $\norm{X}$ is Breiman's. \citet{Breiman1965} proved the
lemma for $0<\beta\leq1$ and \citet{ClineSamorodnitsky1994} for all
$\beta>0$; see \citet{JessenMikosch2006} for a review of products of
regularly varying random variables. For completeness, Appendix~\ref{sec:proofth3} gives the
standard dominated-convergence proof in the present notation, with
Potter's bounds controlling $H$ near zero and the moment condition
controlling large values of $\norm{X}$. The tail exponent $\beta$ is
inherited entirely from the denominator and does not depend on the
dimension $p$, since only $\norm{X}$ enters; the numerator contributes
only the constant $E\norm{X}^\beta$. What the ratio construction adds
are the density-level counterparts at both ends,
Theorems~\ref{thm:origin_power} and~\ref{thm:density_tail} and
Proposition~\ref{prop:origin_log}, which are what the posterior
analysis of Section~\ref{sec:gaussian_sampling} requires.

The most common case is a denominator density with a positive, finite
limit at zero, as for the Cauchy distribution. Then $\beta=1$, and the
moment condition of Theorem~\ref{thm:radial_tail} can be relaxed to
finiteness of the first moment.

\begin{corollary}[Cauchy-order radial tail]
\label{cor:cauchy_tail}
Let $X\in\R^p$ and $Y\in\R$ be independent, with
$0<E\norm{X}<\infty$. Suppose that the density $q$ of $|Y|$ satisfies
$$
q(y)\longrightarrow q_0\in(0,\infty),
\qquad y\downarrow0.
$$
Then
\begin{equation}
P(\norm{X/Y}>r)\sim\frac{q_0E\norm{X}}{r},
\qquad r\to\infty.
\label{eq:cauchy_tail}
\end{equation}
If the density $g$ of $Y$ is continuous at zero, then $q_0=2g(0)$.
\end{corollary}

\begin{proof}
Write $W=\norm{X}$ and $k(t)=H(t)/t$ for $t>0$. By
\eqref{eq:tail_identity},
$$
rP(\norm{X/Y}>r)=E\{Wk(W/r)\}.
$$
Since $H(t)=\int_0^tq(y)\,dy$, we have $k(t)\to q_0$ as $t\downarrow0$.
Moreover, $k$ is bounded on $(0,\infty)$, since it is bounded near zero
and $k(t)\leq1/t$. As $r\to\infty$, $Wk(W/r)\to q_0W$ almost surely,
and the bound $Wk(W/r)\leq W\sup_{t>0}k(t)$, together with $EW<\infty$,
allows dominated convergence:
$rP(\norm{X/Y}>r)\to q_0EW$.
\end{proof}

For the classical horseshoe, $q_0=2/\pi$ and $E|X|=\tau\sqrt{2/\pi}$,
so \eqref{eq:cauchy_tail} gives
$P(|\Theta|>r)\sim2\tau\sqrt{2/\pi}/(\pi r)$, the survival function
corresponding to the density tail \eqref{eq:hs_tail_order} of
Section~\ref{sec:horseshoe_review}.

In most examples the density $q$ is more accessible than $H$, and a
density version of the small-denominator condition is more convenient.
For a general tail index $\beta>0$, suppose that
\begin{equation}
q(y)\sim c_Ty^{\beta-1}L_T(1/y),
\qquad y\downarrow0,
\label{eq:q_tail_rv}
\end{equation}
where $c_T>0$ and $L_T$ is positive and slowly varying. Karamata's
theorem then gives $H(t)\sim(c_T/\beta)t^\beta L_T(1/t)$ as
$t\downarrow0$, so Theorem~\ref{thm:radial_tail} applies with
$c_Y=c_T/\beta$. Under its moment condition, $0<\beta<1$ produces radial
tails heavier than Cauchy order, and $\beta>1$ produces lighter
polynomial tails; in the former case $q$ diverges at zero, and in the
latter it vanishes. Condition \eqref{eq:q_tail_rv} will reappear in
Section~\ref{sec:kappa_geometry}, where it governs the escape endpoint
of the shrinkage coordinate.

A tail statement for the distribution does not by itself determine the
behavior of the density, which is what enters the posterior analysis of
Section~\ref{sec:gaussian_sampling}. The pointwise counterpart of
Theorem~\ref{thm:radial_tail} is as follows.

\begin{theorem}[Directional density tail]
\label{thm:density_tail}
Suppose that \eqref{eq:q_tail_rv} holds and that $q$ is bounded on
$[\delta,\infty)$ for every $\delta>0$. For a unit vector $e\in\R^p$,
assume that, for some $\eta>\max\{0,1-\beta\}$,
$$
\int_0^\infty\left(u^{p+\beta-1-\eta}+u^{p+\beta-1+\eta}\right)
f(ue)\,du<\infty.
$$
Define
$$
C_\beta(e)=\int_0^\infty u^{p+\beta-1}f(ue)\,du.
$$
If $C_\beta(e)>0$, then
\begin{equation}
\pi(re)\sim c_TC_\beta(e)\,r^{-(p+\beta)}L_T(r),
\qquad r\to\infty.
\label{eq:density_tail}
\end{equation}
\end{theorem}

The argument mirrors that of Theorem~\ref{thm:origin_power}, now at the
other end of the denominator: after the substitution $u=ry$ in
\eqref{eq:radial_ratio_density}, the small values of $y$ produce the
factor $r^{-(p+\beta)}L_T(r)$, while boundedness of $q$ away from zero
makes the remaining range negligible. The details are in
Appendix~\ref{sec:proofth4}.

\begin{remark}
The pointwise density has tail exponent $p+\beta$. In polar
coordinates, $d\theta=r^{p-1}\,dr\,d\omega(e)$, so for each direction
$e$ with $C_\beta(e)>0$,
$$
r^{p-1}\pi(re)\sim c_TC_\beta(e)\,r^{-1-\beta}L_T(r),
\qquad r\to\infty,
$$
in agreement with the dimension-free radial survival exponent $\beta$
of Theorem~\ref{thm:radial_tail}.
\end{remark}

The moment condition in Theorem~\ref{thm:radial_tail} is substantive.
When the numerator is itself heavy-tailed, it may contribute a slowly
varying correction to the tail generated by small denominators. The
Cauchy numerator is the leading example: it has $E\norm{X}=\infty$,
so for $\beta=1$ neither Theorem~\ref{thm:radial_tail} nor
Corollary~\ref{cor:cauchy_tail} applies. For $X\in\R^p$, define the
truncated first radial moment
$$
M_X(r)=E\{\norm{X};\,\norm{X}\leq r\}.
$$

\begin{remark}[Borderline heavy numerator]
\label{rem:borderline}
If $H(t)\sim c_Yt$ as $t\downarrow0$, $M_X$ is slowly varying with
$M_X(r)\to\infty$, and $rP(\norm{X}>r)=o\{M_X(r)\}$, then
$P(\norm{X/Y}>r)\sim c_YM_X(r)/r$. The proof splits
$E\{H(\norm{X}/r)\}$ at $\norm{X}=\delta r$ and uses $H(t)\sim c_Yt$ on
the first part, as in Theorem~\ref{thm:radial_tail}. For a Cauchy
numerator, $M_X$ grows logarithmically, which produces the logarithmic
correction in the Cauchy/Cauchy prior of Section~\ref{sec:cauchycauchy}.
\end{remark}

\section{Gaussian Likelihood and Large-Signal Robustness}
\label{sec:gaussian_sampling}

The heavy tails derived in Section~\ref{sec:ratio_tails} are not an end
in themselves: their purpose is to keep large signals from being shrunk
excessively. This section makes that precise under the multivariate
normal means model
\begin{equation}
D\mid\theta\sim N_p(\theta,\sigma^2I_p),
\qquad \sigma>0,
\label{eq:norlik}
\end{equation}
where $\theta\in\R^p$ and $p\geq1$. Write $d$ for an observed value of
$D$, let $\pi$ be a spherically symmetric prior density on $\R^p$, and
denote the posterior mean by
$\widehat\theta_{\rm mean}(d)=E(\Theta\mid D=d)$. Ratio priors with
spherically symmetric numerators (Remark~\ref{rem:mvgeometry}) are the
main example, but the results below use only the radial tail of $\pi$.

The posterior mean is conveniently expressed through the marginal
density
$$
m(d)=\int_{\R^p}\phi_{\sigma,p}(d-\theta)\pi(\theta)\,d\theta,
$$
where $\phi_{\sigma,p}$ denotes the $N_p(0,\sigma^2I_p)$ density.
Differentiation under the integral gives
\begin{equation}
\nabla m(d)=\frac{1}{\sigma^2}\int_{\R^p}
(\theta-d)\phi_{\sigma,p}(d-\theta)\pi(\theta)\,d\theta,
\label{eq:marginal_derivative}
\end{equation}
and consequently Brown's identity \citep{brown1971admissible}, 
\begin{equation}
\widehat\theta_{\rm mean}(d)=d+\sigma^2\nabla\log m(d).
\label{eq:brown_identity}
\end{equation}
For the use of Brown's identity in Bayesian asymptotic
analysis, see \citet[Section~20.9]{DasGupta2008}.
The shrinkage applied to an observation is thus governed entirely by
the score of the marginal, and the large-signal behavior of the
posterior mean reduces to that of $\nabla\log m(d)$ as
$\norm{d}\to\infty$.

To describe the tail, write $\pi(\theta)=h(\norm{\theta})$ and suppose
that, for some $\beta>0$,
\begin{equation}
h(r)=r^{-(p+\beta)}L(r)\{1+o(1)\},
\qquad r\to\infty,
\label{eq:prior_regular_tail_gaussian}
\end{equation}
where $L$ is positive and slowly varying. For ratio priors, this is the
conclusion of Theorem~\ref{thm:density_tail}, with $\beta$ the
small-denominator index in \eqref{eq:q_tail_rv}. Since the score
involves a derivative, we also require a mild smoothness condition:
$h$ is eventually positive and continuously differentiable, with
\begin{equation}
\frac{rh'(r)}{h(r)}\longrightarrow-(p+\beta),
\qquad r\to\infty.
\label{eq:smooth_regular_variation}
\end{equation}
Condition \eqref{eq:smooth_regular_variation} in fact implies
\eqref{eq:prior_regular_tail_gaussian}, with $L$ normalized slowly
varying \citep[Section~1.3]{BinghamGoldieTeugels1987}. If
$h(r)=r^{-(p+\beta)}L(r)$ eventually, it holds whenever $L$ is
eventually differentiable with $rL'(r)/L(r)\to0$. Unlike
\eqref{eq:prior_regular_tail_gaussian}, it is not supplied by the
general theory and must be checked for each prior; for Cauchy/Cauchy
this is done in Section~\ref{sec:cauchycauchy}.

Under these conditions, Gaussian smoothing does not alter the tail of
the prior, either in size or in score.

\begin{lemma}[Gaussian convolution in the tail]
\label{lem:gaussian_convolution_tail}
Let $\pi$ be a spherically symmetric probability density satisfying
\eqref{eq:prior_regular_tail_gaussian} and
\eqref{eq:smooth_regular_variation} for some $\beta>0$. Then, as
$\norm{d}\to\infty$,
\begin{equation}
m(d)\sim\pi(d),
\label{eq:marginal_prior_equivalence}
\end{equation}
and
\begin{equation}
\nabla\log m(d)=-(p+\beta)\frac{d}{\norm{d}^2}+o(\norm{d}^{-1}).
\label{eq:marginal_score_tail}
\end{equation}
\end{lemma}

The idea is that, for large $\norm{d}$, the Gaussian kernel averages
$\pi$ over a neighborhood on which a regularly varying density is
nearly constant; the smoothness condition controls the first-order
variation that produces the score. The proof is given in
Appendix~\ref{sec:prooflemma}. Combined with \eqref{eq:brown_identity},
the lemma yields the main result of this section.

\begin{theorem}[Large-signal robustness]
\label{thm:posterior_mean_robustness}
Under the assumptions of Lemma~\ref{lem:gaussian_convolution_tail},
the posterior mean satisfies
\begin{equation}
\widehat\theta_{\rm mean}(d)
=d-\frac{\sigma^2(p+\beta)}{\norm{d}^2}\,d+o(\norm{d}^{-1}),
\qquad \norm{d}\to\infty.
\label{eq:posterior_mean_expansion}
\end{equation}
In particular,
$$
\norm{\widehat\theta_{\rm mean}(d)-d}\longrightarrow0,
\qquad
\frac{\norm{\widehat\theta_{\rm mean}(d)}}{\norm{d}}\longrightarrow1.
$$
Under squared-error loss, the corresponding frequentist risk satisfies
\begin{equation}
R(\theta,\widehat\theta_{\rm mean})
=E_\theta\left[\norm{\widehat\theta_{\rm mean}(D)-\theta}^2\right]
\longrightarrow p\sigma^2,
\qquad \norm{\theta}\to\infty.
\label{eq:large_signal_risk}
\end{equation}
\end{theorem}

The expansion \eqref{eq:posterior_mean_expansion} is immediate from
\eqref{eq:brown_identity} and \eqref{eq:marginal_score_tail}. The risk
statement follows because the correction
$\widehat\theta_{\rm mean}(d)-d$ is bounded and vanishes at infinity;
the details are in Appendix~\ref{sec:proofth5}.

The theorem says that large signals are left essentially untouched:
the correction vanishes at rate $\norm{d}^{-1}$, and the risk
approaches $p\sigma^2$, that of the unshrunk observation $D$. The size
of the correction, $\sigma^2(p+\beta)/\norm{d}$, decreases with the
tail weight of the prior, so heavier tails, generated by small
denominators, leave large signals less attenuated. For a Laplace prior
with scale $b$, by contrast, the
posterior mean retains a bias of size approaching $\sigma^2/b$, and the
scalar risk tends to $\sigma^2+\sigma^4/b^2$ rather than $\sigma^2$.

\section{Shrinkage Geometry and the Horseshoe Link}
\label{sec:kappa_geometry}

Sections~\ref{sec:ratio_priors}--\ref{sec:gaussian_sampling} concern the
marginal prior: its spike, its tails, and their consequences for large
signals. The horseshoe is also known for a third feature, the U-shaped
distribution of its shrinkage coefficient. We now show that the ratio
construction produces this feature as well, and that it is generated by
the denominator alone.

For $\Theta=X/Y$, large values of $|Y|$ contract the ratio toward zero,
whereas values near zero allow large ratios. The bounded coordinate
\begin{equation}
\kappa_Y=\frac{Y^2}{1+Y^2}
\label{eq:kappaY_definition}
\end{equation}
maps escape, $|Y|\downarrow0$, to $\kappa_Y\downarrow0$, and
contraction, $|Y|\to\infty$, to $\kappa_Y\uparrow1$. Because $\kappa_Y$
depends only on the denominator, its distribution describes the
contraction--escape geometry of the ratio representation independently
of the numerator and of any likelihood; its relationship to posterior
shrinkage is examined in Section~\ref{sec:kappa_posterior}. With $q$ the
density of $|Y|$, a change of variables gives the density of $\kappa_Y$,
\begin{equation}
h_Y(\kappa)=
\frac{q\left\{\left(\kappa/(1-\kappa)\right)^{1/2}\right\}}
{2\kappa^{1/2}(1-\kappa)^{3/2}},
\qquad 0<\kappa<1,
\label{eq:kappaY_density}
\end{equation}
so the behavior of $q$ near zero determines $h_Y$ near $\kappa=0$,
while the tail of $q$ determines $h_Y$ near $\kappa=1$. The two indices
of Sections~\ref{sec:ratio_priors} and~\ref{sec:ratio_tails} therefore
reappear, now as the exponents of the two arms of $h_Y$.

\begin{theorem}[Denominator-induced horseshoe geometry]
\label{thm:kappa_endpoint_geometry}
Suppose that \eqref{eq:q_tail_rv} and \eqref{eq:q_origin_rv} hold with
indices $\beta>0$ and $\alpha>0$, respectively. Then
$$
\begin{aligned}
h_Y(\kappa)&\sim\frac{c_T}{2}\,\kappa^{\beta/2-1}L_T(\kappa^{-1/2}),
&&\kappa\downarrow0,\\
h_Y(\kappa)&\sim\frac{c_O}{2}\,(1-\kappa)^{\alpha/2-1}
L_O\{(1-\kappa)^{-1/2}\},
&&\kappa\uparrow1.
\end{aligned}
$$
Consequently, $0<\beta<2$ implies $h_Y(\kappa)\to\infty$ at the escape
endpoint $\kappa=0$, while $0<\alpha<2$ implies $h_Y(\kappa)\to\infty$
at the contraction endpoint $\kappa=1$.
\end{theorem}

\begin{proof}
Substitute $y=\{\kappa/(1-\kappa)\}^{1/2}$ in \eqref{eq:kappaY_density}
and use $y\sim\kappa^{1/2}$ as $\kappa\downarrow0$ and
$y\sim(1-\kappa)^{-1/2}$ as $\kappa\uparrow1$.
\end{proof}

The tail index $\beta$ thus controls the escape arm, and the origin
index $\alpha$ the contraction arm. For $\alpha=\beta=1$ the endpoint
exponents are those of the ${\rm Beta}(1/2,1/2)$ horseshoe density, and
the correspondence becomes exact for the family
$$
q_{a,b}(y)=\frac{2}{B(a,b)}\frac{y^{2a-1}}{(1+y^2)^{a+b}},
\qquad y>0,
$$
for which $\kappa_Y\sim{\rm Beta}(a,b)$, with $\beta=2a$ and
$\alpha=2b$. The standard Cauchy denominator corresponds to $a=b=1/2$:
\begin{equation}
Y\sim C(0,1)\quad\Longrightarrow\quad
\kappa_Y\sim{\rm Beta}\left(\frac12,\frac12\right).
\label{eq:cauchy_kappa_exact}
\end{equation}
The characteristic U-shaped density is therefore generated by the
denominator alone and does not require a Gaussian numerator.

\begin{figure}[ht]
\centering
\begin{minipage}[t]{0.48\textwidth}
\centering
\includegraphics[width=\textwidth]
{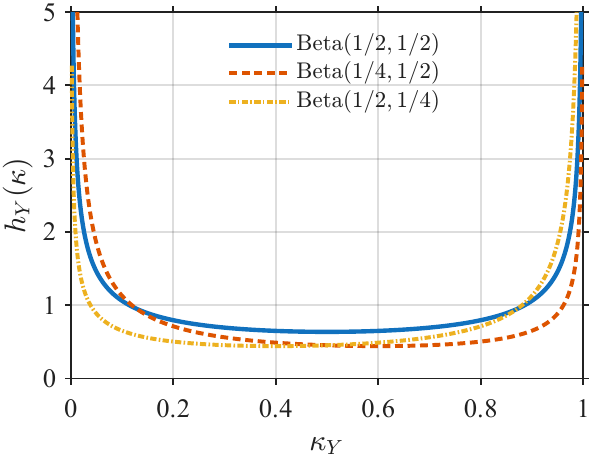}

{\small (a)}
\end{minipage}
\hfill
\begin{minipage}[t]{0.48\textwidth}
\centering
\includegraphics[width=\textwidth]
{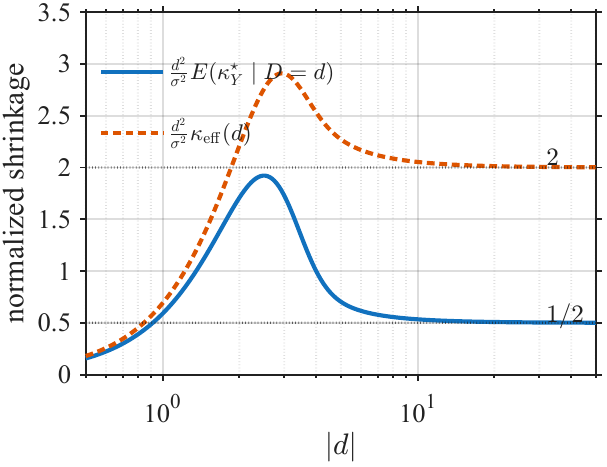}

{\small (b)}
\end{minipage}
\caption{Shrinkage geometry of ratio priors.
(a) Denominator-induced Beta densities for
$(\alpha,\beta)=(1,1),(1,1/2),(1/2,1)$.
(b) For Uniform/Cauchy, the normalized latent shrinkage
$d^2E(\kappa_Y^\star\mid d)/\sigma^2$ and effective shrinkage
$d^2\kappa_{\rm eff}(d)/\sigma^2$ approach $1/2$ and $2$,
respectively ($\tau=\sigma=1$).}
\label{fig:kappa_geometry}
\end{figure}

Figure~\ref{fig:kappa_geometry}a shows how the two arms can be varied
separately. The classical case $(\alpha,\beta)=(1,1)$ has equal arms.
Keeping $\alpha=1$ but reducing $\beta$ to $1/2$ strengthens the
$\kappa=0$ arm and, under the moment condition of
Theorem~\ref{thm:radial_tail}, produces a tail heavier than Cauchy
order. Keeping $\beta=1$ but reducing $\alpha$ to $1/2$ strengthens the
$\kappa=1$ arm and, in the scalar case, replaces the logarithmic origin
singularity by a power singularity. The two features of horseshoe
geometry can therefore be adjusted independently.

\subsection{From denominator geometry to posterior shrinkage}
\label{sec:kappa_posterior}

The coordinate $\kappa_Y$ is defined without reference to data. To
relate it to posterior shrinkage, we calibrate it to the sampling
model. For a symmetric numerator write $X=\tau Z$ with $\tau>0$ and,
under the Gaussian sampling model \eqref{eq:norlik}, define the
scale-calibrated coordinate
\begin{equation}
\kappa_Y^\star=\frac{\sigma^2Y^2}{\tau^2+\sigma^2Y^2}.
\label{eq:kappaY_star}
\end{equation}
It is a monotone reparameterization of $\kappa_Y$ with the same
endpoints and endpoint exponents, and $\kappa_Y^\star=\kappa_Y$ when
$\sigma=\tau$. Conditionally on the denominator,
$$
\Theta\mid\kappa_Y^\star
\stackrel{d}{=}
\sigma\left(\frac{1-\kappa_Y^\star}{\kappa_Y^\star}\right)^{1/2}Z,
$$
so $\kappa_Y^\star$ is an exact latent contraction scale for every
numerator: $\kappa_Y^\star\uparrow1$ collapses the conditional prior
toward zero, whereas $\kappa_Y^\star\downarrow0$ makes it diffuse.

When $Z\sim N_p(0,I_p)$, the conditional prior is
$N_p(0,\tau^2\lambda^2I_p)$ with $\lambda=1/|Y|$, and
$\kappa_Y^\star=\sigma^2/(\sigma^2+\tau^2\lambda^2)$ is the coefficient
in \eqref{eq:intro_shrinkage_factor}. Normal conjugacy gives
\begin{equation}
\widehat\theta_{\rm mean}(d)=\{1-E(\kappa_Y^\star\mid D=d)\}\,d.
\label{eq:kappa_gaussian_identity}
\end{equation}
If, moreover, $\kappa_Y\sim{\rm Beta}(a,b)$, a change of variables
shows that $\kappa_Y^\star$ has density
$$
\frac{\varphi^{a}}{B(a,b)}\,
\kappa^{a-1}(1-\kappa)^{b-1}\{1+(\varphi-1)\kappa\}^{-(a+b)},
\qquad 0<\kappa<1,\quad \varphi=\tau^2/\sigma^2,
$$
the three-parameter beta distribution of \citet{armagan2011generalized},
with $a$ and $b$ interchanged relative to their notation. For $p=1$, a
standard Cauchy denominator makes $\lambda$ half-Cauchy; when
$\sigma=\tau$, the resulting conditional model is the classical
horseshoe, and $\kappa_Y^\star\sim{\rm Beta}(1/2,1/2)$.

What is special about the Gaussian numerator is not the latent geometry
but the linear identity \eqref{eq:kappa_gaussian_identity}. To compare
latent and actual shrinkage for general numerators, note that for a
spherically symmetric prior the posterior mean is radial,
$\widehat\theta_{\rm mean}(d)=a_p(\norm{d})\,d$, and define the
effective posterior-mean attenuation
\begin{equation}
\kappa_{\rm eff}(d)=1-a_p(\norm{d}),\qquad d\ne0,
\label{eq:kappa_effective}
\end{equation}
which equals $1-\widehat\theta_{\rm mean}(d)/d$ when $p=1$. By
\eqref{eq:kappa_gaussian_identity},
$\kappa_{\rm eff}(d)=E(\kappa_Y^\star\mid d)$ for a Gaussian numerator;
for non-Gaussian numerators the latent contraction and the effective
shrinkage generally differ.

\subsection{Large-signal geometry}

Theorem~\ref{thm:posterior_mean_robustness} already fixes the behavior
of the effective shrinkage: $\kappa_{\rm eff}(d)$ is of order $d^{-2}$.
It remains to compare the latent coordinate with it, and the comparison
is sharpest for large signals.

\begin{theorem}[Large-signal latent and effective shrinkage]
\label{thm:kappa_large_signal}
Let $p=1$ and $X=\tau Z$, where $Z$ has a symmetric density $f_Z$.
Suppose that \eqref{eq:q_tail_rv} holds, that $q$ is bounded on
$[\delta,\infty)$ for every $\delta>0$, and that
$E|Z|^{\beta+2+\eta}<\infty$ for some $\eta>0$. Let $\mu_r=E|Z|^r$ and
let $Z_\beta$ have the $\beta$-power-biased density
$f_\beta(z)=|z|^\beta f_Z(z)/\mu_\beta$. Then, as $|\theta|\to\infty$,
$$
\theta^2\kappa_Y^\star\mid\{\Theta=\theta\}
\Longrightarrow\sigma^2Z_\beta^2,
\qquad
E(\kappa_Y^\star\mid\Theta=\theta)\sim
\frac{\sigma^2}{\theta^2}\frac{\mu_{\beta+2}}{\mu_\beta},
$$
and, as $|d|\to\infty$,
\begin{equation}
E(\kappa_Y^\star\mid D=d)\sim
\frac{\sigma^2}{d^2}\frac{\mu_{\beta+2}}{\mu_\beta}.
\label{eq:kappa_latent_large}
\end{equation}
If, in addition, the induced density $\pi$ satisfies the smoothness
condition \eqref{eq:smooth_regular_variation} with $p=1$, then
\begin{equation}
\kappa_{\rm eff}(d)=\frac{\sigma^2(1+\beta)}{d^2}+o(d^{-2})
\label{eq:kappa_eff_large}
\end{equation}
and
\begin{equation}
\frac{E(\kappa_Y^\star\mid D=d)}{\kappa_{\rm eff}(d)}
\longrightarrow
\frac{\mu_{\beta+2}}{(1+\beta)\mu_\beta}.
\label{eq:kappa_ratio_large}
\end{equation}
\end{theorem}

The mechanism is a tilting of the numerator. A large value of $\Theta$
is explained mainly by a small denominator, and conditionally on
$\Theta=\theta$ the numerator variable $|Z|$ acquires the
$\beta$-power-biased law; the moment ratio in \eqref{eq:kappa_ratio_large}
records this tilt. The proof is given in Appendix~\ref{sec:large_signal}.

For $Z\sim N(0,1)$, $\mu_{\beta+2}=(1+\beta)\mu_\beta$, so the limit in
\eqref{eq:kappa_ratio_large} equals one, consistently with
\eqref{eq:kappa_gaussian_identity}. For non-Gaussian numerators, the
two coefficients share the quadratic rate but generally have different
constants.

\begin{example}[Uniform/Cauchy]
\label{ex:unif_cauchy}
Let $Z\sim{\rm Unif}(-1,1)$ and $Y\sim C(0,1)$, as in the third row of
Table~\ref{tab:ratio_examples}. Then $\beta=1$ and $\mu_3/\mu_1=1/2$, so
$$
E(\kappa_Y^\star\mid d)\sim\frac{\sigma^2}{2d^2},
\qquad
\kappa_{\rm eff}(d)\sim\frac{2\sigma^2}{d^2}.
$$
Figure~\ref{fig:kappa_geometry}b shows this separation: the two
quantities share the $d^{-2}$ order, but their constants retain
information about the numerator. The induced density is
$\pi(\theta)=(2\pi\tau)^{-1}\log(1+\tau^2/\theta^2)$, the horseshoe-type
prior of \citet{Bhadra2021}, which arises in Section~\ref{sec:cauchycauchy}
from a Cauchy numerator and a Pareto denominator. The same prior thus
has two ratio representations with different latent coordinates: here
$\kappa_Y\sim{\rm Beta}(1/2,1/2)$, whereas under the Pareto denominator
$|Y|>1$ and $\kappa_Y$ is confined to $(1/2,1)$. Latent shrinkage
geometry therefore depends on the particular ratio representation,
whereas effective posterior shrinkage is determined by the induced prior
and the likelihood.
\end{example}

The Cauchy/Cauchy construction of Section~\ref{sec:cauchycauchy} provides
the complementary case, in which the numerator moment condition of
Theorem~\ref{thm:kappa_large_signal} fails; there the latent and
effective coefficients differ not only in their constants but in their
rates.

\section{Representative Ratio Priors and the Cauchy/Cauchy Example}
\label{sec:cauchycauchy}

The theory of the preceding sections applies to a broad range of ratio
constructions. Table~\ref{tab:ratio_examples} summarizes representative
examples and shows how the numerator and denominator jointly determine
concentration near zero and tail behavior. We then examine the
Cauchy/Cauchy prior in detail. It brings together several strands of
the general theory: dimension-dependent origin singularities,
logarithmically modified Cauchy-order tails, large-signal robustness,
and the distinction between denominator-induced geometry and effective
posterior shrinkage.

\begin{table}[htbp]
\centering
\small
\renewcommand{\arraystretch}{1.25}
\begin{tabular}{@{}
>{\raggedright\arraybackslash}p{0.25\textwidth}
>{\raggedright\arraybackslash}p{0.23\textwidth}
>{\raggedright\arraybackslash}p{0.21\textwidth}
>{\raggedright\arraybackslash}p{0.21\textwidth}@{}}
\hline
Numerator $X$ & Denominator $Y$ & Behavior at the origin & Tail behavior \\ \hline
$N(0,\tau^2)$
& $C^+(0,1)$
& logarithmic singularity
& Cauchy order; classical horseshoe \\
$N_p(0,\tau^2\Sigma)$
& $C(0,1)$
& logarithmic for $p=1$; $\pi(re)\asymp r^{1-p}$ for $p>1$
& radial survival of order $r^{-1}$; shared-scale block prior \\
${\rm Unif}(-\tau,\tau)$
& $C(0,1)$
& logarithmic singularity
& Cauchy order; prior of \citet{Bhadra2021} \\
$N(0,\tau^2)$
& Double Pareto$(\alpha)$
& power singularity for $0<\alpha<1$; logarithmic at $\alpha=1$; finite for $\alpha>1$
& Cauchy-order tail \\
$N(0,\tau^2)$
& $1/\lambda$,~ $\lambda^2\sim\beta'(a,b)$
& power singularity for $a<1/2$; logarithmic at $a=1/2$; finite for $a>1/2$
& survival order $r^{-2b}$; normal--beta prime prior of \citet{BaiGhosh2019, BaiGhosh2021} \\
$t_\nu(0,\tau)$
& $t_\nu(0,1)$, $0<\nu\leq1$
& power singularity for $0<\nu<1$; logarithmic at $\nu=1$
& survival order $r^{-\nu}$ for $0<\nu<1$; $(\log r)/r$ for $\nu=1$ \\
$\tau S_1U^{1/c}$, $U\sim{\rm Ga}(a,1)$
& $S_2V^{1/c}$, $V\sim{\rm Ga}(b,1)$
& numerator-generated singularity if $ac<1$
& Cauchy order when $bc=1$ \\
$C(0,\tau)$
& $C(0,1)$
& logarithmic singularity
& density $\asymp(\log|\theta|)/\theta^2$ \\ \hline
\end{tabular}
\caption{Representative ratio-generated priors. Entries refer to the
scalar case unless stated otherwise. In the normal--beta prime row,
$\beta'(a,b)$ denotes the beta prime distribution with density
$x^{a-1}(1+x)^{-(a+b)}/B(a,b)$, $x>0$. In the signed Gamma-power row,
$a,b,c>0$, and $S_1,S_2$ are independent Rademacher variables,
independent of $U,V$. The Double Pareto density is
$g_\alpha(y)=\frac12(1+|y|/\alpha)^{-(\alpha+1)}$. Cauchy order means
radial survival proportional to $r^{-1}$, up to any indicated slowly
varying factor.}
\label{tab:ratio_examples}
\end{table}

The first row recovers the classical horseshoe conditional on the
global scale $\tau$, since the reciprocal of a standard half-Cauchy
variable is again half-Cauchy. The multivariate normal/Cauchy
construction introduces a common scale across coordinates and thus a
mechanism for joint shrinkage; in the spherical case the resulting
posterior attenuation is radial. The Double Pareto family controls the
origin singularity through the denominator tail index; for the
normal--beta prime row, $\alpha=2a$ and $\beta=2b$, so the spike
condition $a\leq1/2$ of \citet{BaiGhosh2019} is exactly $\alpha\leq1$. The Student-$t$ example shows how a heavy-tailed
numerator can determine the tail order, and in the signed Gamma-power
construction the origin singularity arises from the numerator instead.

Ratio representations need not be unique. The horseshoe-type prior of
\citet{Bhadra2021}, induced by the uniform/Cauchy row
(Example~\ref{ex:unif_cauchy}), also arises from
a Cauchy numerator $X\sim C(0,\tau)$ and an independent symmetric Pareto
Type~I denominator with $g(y)=(2y^2)^{-1}\mathbf{1}(|y|>1)$. In both
representations the logarithmic origin singularity is the boundary case
$\alpha=p=1$ of Proposition~\ref{prop:origin_log}. The tail, however, is
attributed to different mechanisms: in the uniform/Cauchy representation
it arises from small denominators (Corollary~\ref{cor:cauchy_tail}),
whereas the Pareto denominator is bounded away from zero, so the tail is
inherited from the numerator.

A denominator with only logarithmic tails lies outside this framework
but marks its natural boundary, as the following remark records.

\begin{remark}[The logarithmic boundary $\alpha=\beta=0$]
\label{rem:logcauchy}
Let $X\sim N(0,\tau^2)$ and $\log|Y|\sim C(0,1)$ independently. Then
$q(y)=\{\pi y(1+\log^2y)\}^{-1}$ is regularly varying with index $-1$ at
both ends, so $\alpha=\beta=0$, outside the positive-index assumptions
of Sections~\ref{sec:ratio_priors}--\ref{sec:ratio_tails}. Several
results nevertheless extend. Theorem~\ref{thm:origin_singularity}
applies, since $E|Y|=\infty$. The proof of Theorem~\ref{thm:origin_power}
uses only $\alpha<p$ and therefore covers $\alpha=0$; with
$L_O(y)=\{\pi\log^2y\}^{-1}$ and $C_0=\int_0^\infty f=1/2$ it gives the
first line below. The second follows by a direct calculation, since
Theorem~\ref{thm:density_tail} requires $\beta>0$:
$$
\pi(\theta)\sim
\frac{1}{2\pi|\theta|}
\begin{cases}
\log^{-2}(\tau/|\theta|), & \theta\to0,\\
\log^{-2}(|\theta|/\tau), & |\theta|\to\infty.
\end{cases}
$$
This is a proper prior with a logarithmic correction to the
scale-invariant density $1/|\theta|$; absolute moments of every nonzero
order are infinite. Since $H(t)\sim\{\pi\log(1/t)\}^{-1}$ is slowly
varying, Breiman's lemma, whose original statement covers index zero
\citep{Breiman1965}, gives $P(|\Theta|>r)\sim\{\pi\log(r/\tau)\}^{-1}$.
Equivalently, concentration under the global scale is only logarithmic:
$P(|\Theta|>c)\sim1/\{\pi\log(c/\tau)\}$ as $\tau\downarrow0$, for fixed
$c>0$. The denominator-induced coordinate has density
$$
h_Y(\kappa)=
\frac{1}{
2\pi\kappa(1-\kappa)
\{1+\frac14\log^2(\kappa/(1-\kappa))\}},
$$
a proper analogue of the Haldane prior. The proof of
Lemma~\ref{lem:gaussian_convolution_tail} also applies
with $p+\beta=1$, since differentiation under the integral sign gives
$rh'(r)/h(r)\to-1$. Hence
$$
\widehat\theta_{\rm mean}(d)
=d-\frac{\sigma^2}{d}+o(|d|^{-1}),
\qquad |d|\to\infty.
$$
The leading correction is half that of the classical horseshoe.
\end{remark}

\subsection{The spherical Cauchy/Cauchy prior}
\label{sec:caucau_prior}

Let $X=\tau Z$, where $Z$ has the standard spherical $p$-variate Cauchy
distribution, and let $Y\sim C(0,1)$ be independent of $X$. The
numerator density is
$$
f(x)=\frac{c_p\tau}{(\tau^2+\norm{x}^2)^{(p+1)/2}},
\qquad
c_p=\frac{\Gamma((p+1)/2)}{\pi^{(p+1)/2}},\qquad p\geq1.
$$
The standard Cauchy denominator has $\alpha=\beta=1$. Its far tail
produces the origin singularity of Section~\ref{sec:ratio_priors}. Its
density has a positive limit at zero, which would give a radial tail of
order $r^{-1}$ under Theorem~\ref{thm:radial_tail}; but
$E\norm{X}=\infty$, so the moment condition of that theorem fails, and
the prior falls in the borderline case of Remark~\ref{rem:borderline}.

Substitution into \eqref{eq:ratio_density} gives, for
$r=\norm{\theta}>0$,
\begin{equation}
\begin{aligned}
\pi_{p,\tau}(\theta)
&=\frac{2c_p}{\pi\tau^p}
\int_0^\infty
\frac{y^p\,dy}{(1+y^2)(1+y^2r^2/\tau^2)^{(p+1)/2}}\\
&=\frac{2c_p}{\pi\tau^p(p+1)}
\,{}_2F_1\!\left(\frac{p+1}{2},\frac{p+1}{2};
\frac{p+3}{2};1-\frac{r^2}{\tau^2}\right).
\end{aligned}
\label{eq:caucau_mv_density}
\end{equation}
For $p=1$ the integral can be evaluated explicitly:
\begin{equation}
\pi_{1,\tau}(\theta)=
\frac{2\tau\log(|\theta|/\tau)}{\pi^2(\theta^2-\tau^2)},
\qquad \theta\ne0,\quad |\theta|\ne\tau,
\label{eq:caucau_density}
\end{equation}
with removable singularities at $\theta=\pm\tau$, where
$\pi_{1,\tau}(\pm\tau)=1/(\pi^2\tau)$. These expressions make the
general results of Sections~\ref{sec:ratio_priors}
and~\ref{sec:ratio_tails} fully explicit, at the origin and in the tail.

Near the origin, with $c_O=2/\pi$, Proposition~\ref{prop:origin_log} for $p=1$ and
Theorem~\ref{thm:origin_power} for $p>1$ give
\begin{equation}
\pi_{p,\tau}(\theta)\sim
\begin{cases}
\displaystyle \frac{1}{\pi^2\tau}
\log\!\left(\frac{\tau^2}{\theta^2}\right), & p=1,\\[7pt]
\displaystyle \frac{2c_p}{\pi\tau(p-1)}
\norm{\theta}^{1-p}, & p>1,
\end{cases}
\qquad \norm{\theta}\to0.
\label{eq:caucau_mv_origin}
\end{equation}
The singularity is logarithmic for $p=1$ and of power order for $p>1$,
and it is integrable in every dimension; for $p=2$ the leading term is
$1/(\pi^2\tau\norm{\theta})$. Figure~\ref{fig:caucau_mv}(a) shows the
bivariate density and its pronounced singularity at the origin.

For the tail, substituting $u=yr/\tau$ in \eqref{eq:caucau_mv_density} gives, for a
unit vector $e$,
\begin{equation}
\pi_{p,\tau}(re)=\frac{2c_p\tau}{\pi r^{p+1}}\,J_p(r/\tau),
\qquad
J_p(a)=\int_0^\infty
\frac{u^p\,du}{(1+u^2)^{(p+1)/2}(1+u^2/a^2)}.
\label{eq:caucau_J_representation}
\end{equation}
Since $u^p(1+u^2)^{-(p+1)/2}\sim u^{-1}$ as $u\to\infty$,
$J_p(a)\sim\log a$, so
\begin{equation}
\pi_{p,\tau}(\theta)\sim
\frac{2c_p\tau}{\pi}
\frac{\log(\norm{\theta}/\tau)}{\norm{\theta}^{p+1}},
\qquad \norm{\theta}\to\infty,
\label{eq:caucau_tail}
\end{equation}
and radial integration gives
\begin{equation}
P(\norm{\Theta}>r)\sim
\frac{2c_p\tau\omega_{p-1}}{\pi}
\frac{\log(r/\tau)}{r},
\qquad r\to\infty,
\label{eq:caucau_mv_survival}
\end{equation}
where $\omega_{p-1}=2\pi^{p/2}/\Gamma(p/2)$ is the surface area of the
unit sphere. The radial survival exponent is one in every dimension;
the logarithmic factor is the borderline contribution of the numerator.
Indeed, $H(t)\sim2t/\pi$ and
$M_X(r)\sim c_p\tau\omega_{p-1}\log(r/\tau)$, and
Remark~\ref{rem:borderline} reproduces
\eqref{eq:caucau_mv_survival}, constant included.

The same representation also verifies the regularity condition required for
posterior robustness. Differentiating $J_p$ under the integral sign and
using the substitution $u=av$ gives $aJ_p'(a)\to1$. With
$h(r)=\pi_{p,\tau}(re)$,
$$
\frac{rh'(r)}{h(r)}
=-(p+1)+\frac{aJ_p'(a)}{J_p(a)}\longrightarrow-(p+1),
$$
so the smooth regular-variation condition
\eqref{eq:smooth_regular_variation} holds with $\beta=1$, and
Theorem~\ref{thm:posterior_mean_robustness} applies.

\subsection{Posterior shrinkage and latent contraction}
\label{sec:caucau_posterior}

Under the Gaussian model \eqref{eq:norlik}, consider the posterior mean
$\widehat\theta_{\rm mean}(d)=E(\Theta\mid D=d)$ and the posterior
spatial median
$\widehat\theta_{\rm med}(d)=\arg\min_{a\in\R^p}E(\norm{\Theta-a}\mid D=d)$,
which reduces to the ordinary posterior median when $p=1$. By rotational
symmetry both estimators are radial:
$$
\widehat\theta_{\rm mean}(d)=a_p(\norm{d})\,d,
\qquad
\widehat\theta_{\rm med}(d)=b_p(\norm{d})\,d,
\qquad d\ne0.
$$
Because the denominator is shared, the attenuation applied to each
coordinate depends on $\norm{d}$ rather than on that coordinate alone.

With $\beta=1$, Theorem~\ref{thm:posterior_mean_robustness} gives
\begin{equation}
\widehat\theta_{\rm mean}(d)
=d-\frac{(p+1)\sigma^2}{\norm{d}^2}\,d
+o(\norm{d}^{-1}),
\qquad \norm{d}\to\infty,
\label{eq:caucau_mv_mean_tail}
\end{equation}
so the posterior mean approaches the observation for large signals, and
its squared-error risk approaches $p\sigma^2$.

\begin{figure}[t]
\centering
\includegraphics[width=0.56\textwidth,
 trim=1cm 1cm 4cm 6cm,clip]{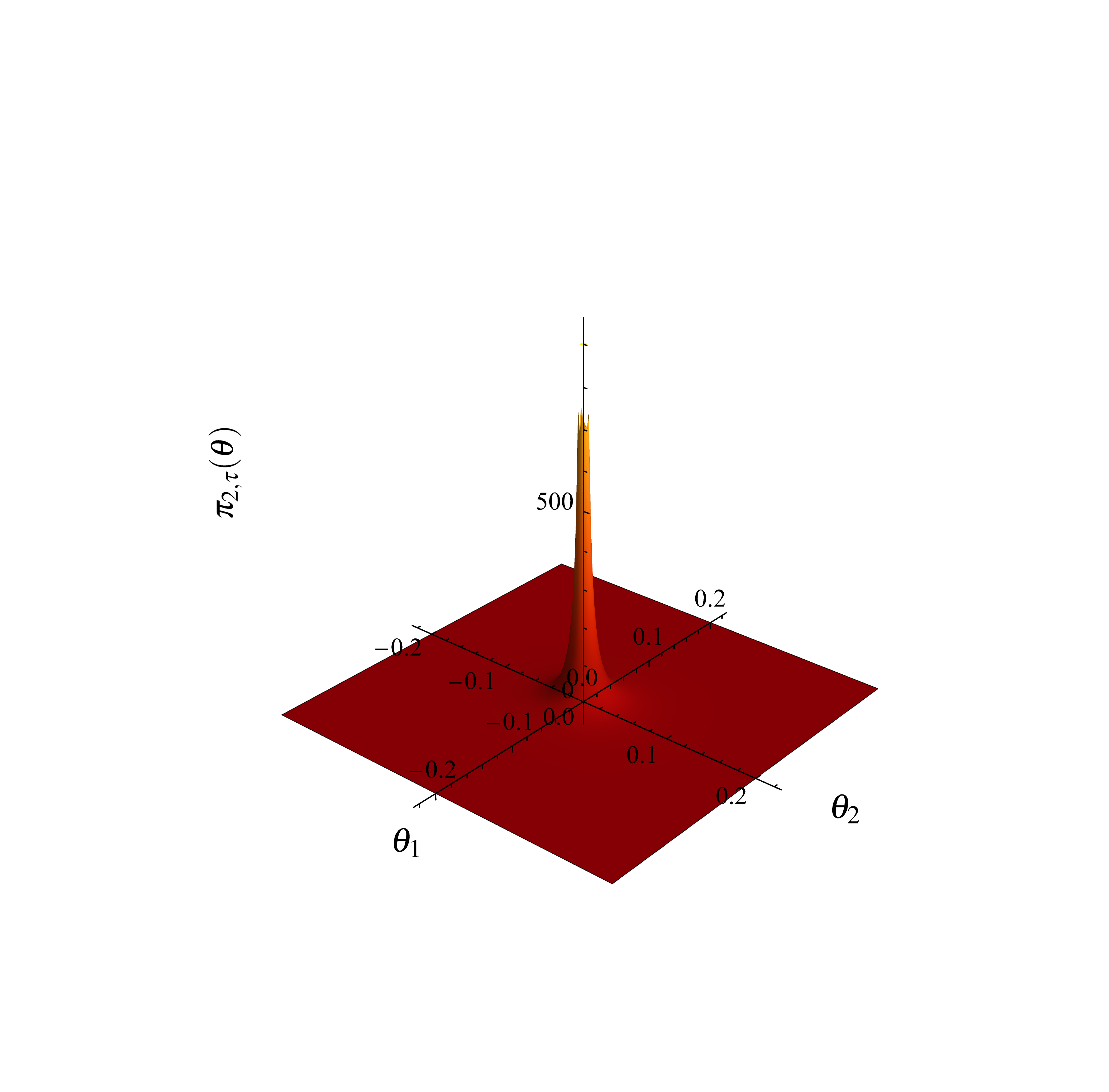}
\hfill
\includegraphics[width=0.42\textwidth]{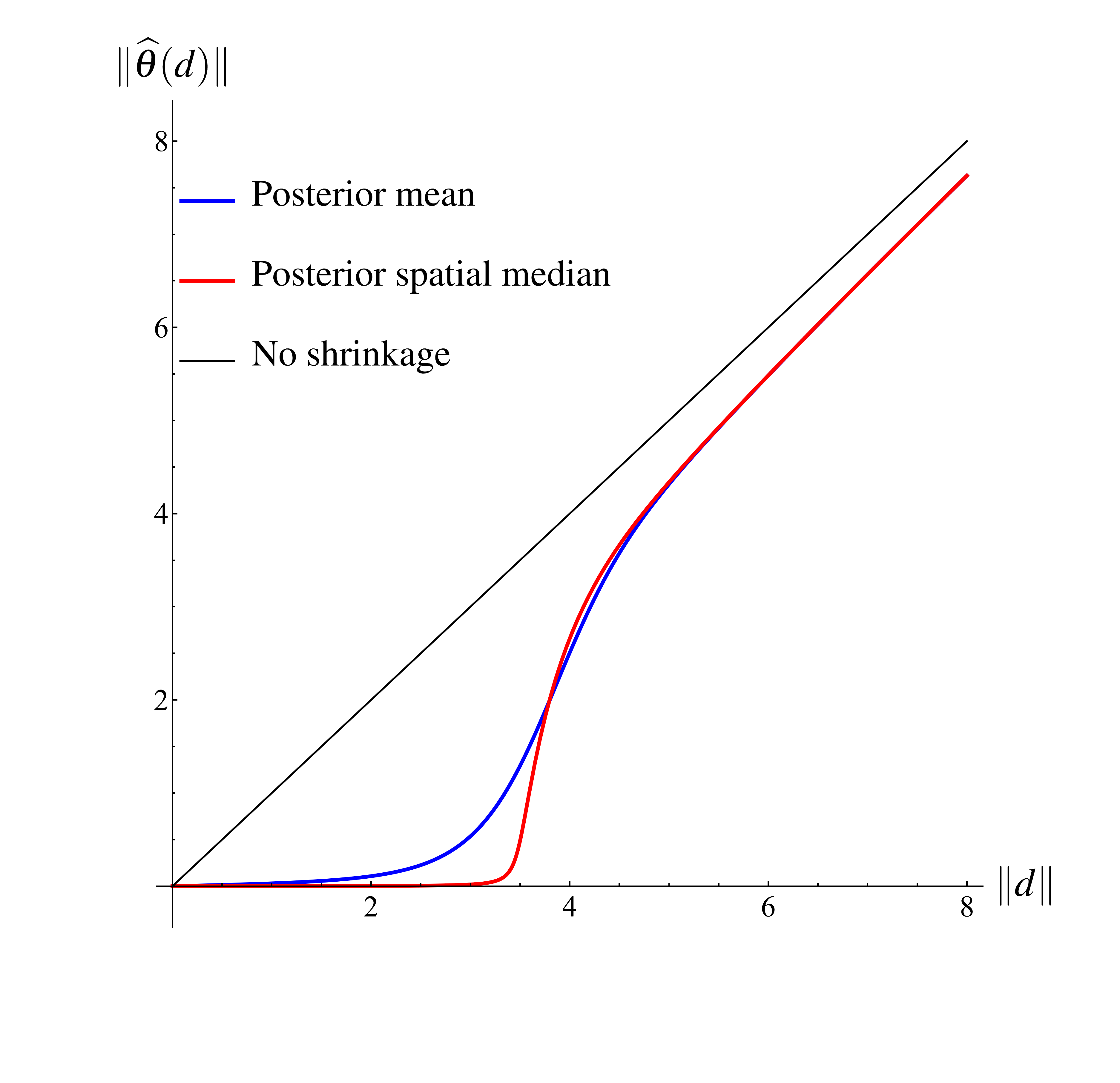}
\caption{Bivariate Cauchy/Cauchy prior and posterior estimators, for
$\tau =1/50$ and $\sigma^2=1$.
(a) Prior density, showing its singularity at the origin.
(b) Magnitudes of the posterior mean and spatial median versus
$\norm{d}$, compared with the identity rule.}
\label{fig:caucau_mv}
\end{figure}

Figure~\ref{fig:caucau_mv}(b) compares the posterior mean and spatial
median for $p=2$, $\tau=1/50$, and $\sigma^2=1$. The posterior mean
attenuates smoothly, whereas the spatial median remains close to zero
over a substantial range of observations before rising sharply,
behaving nearly as a threshold rule. The asymptotic result above
concerns the posterior mean; we do not pursue the corresponding median
asymptotics here.

The shrinkage geometry of Section~\ref{sec:kappa_geometry} adds a
further layer. By \eqref{eq:cauchy_kappa_exact}, the Cauchy denominator
gives $\kappa_Y\sim{\rm Beta}(1/2,1/2)$ in every dimension. For a Cauchy
numerator, however, $\kappa_Y$ is not the Gaussian mixing coefficient.
Writing
$$
X=\frac{\tau Z_0}{\sqrt V},\qquad
Z_0\sim N_p(0,I_p),\quad V\sim\chi_1^2,
$$
with $Z_0$, $V$, and $Y$ mutually independent, gives
$\Theta\mid V,Y\sim N_p\{0,\tau^2(VY^2)^{-1}I_p\}$. The Gaussian mixing
coefficient is therefore
\begin{equation}
\kappa_{VY}=\frac{\sigma^2VY^2}{\tau^2+\sigma^2VY^2},
\qquad
\widehat\theta_{\rm mean}(d)
=\{1-E(\kappa_{VY}\mid D=d)\}\,d.
\label{eq:caucau_kappa_star}
\end{equation}
Both $\kappa_Y$ and its calibrated version $\kappa_Y^\star$ in
\eqref{eq:kappaY_star} omit the numerator's mixing variable $V$.

The difference persists for large observations. With
$a=\norm{\theta}/\tau$, the conditional law of $|Y|$ given
$\Theta=\theta$ yields the exact identity
$$
E(\kappa_Y\mid\Theta=\theta)=\frac{aJ_p'(a)}{2J_p(a)},
$$
so $E(\kappa_Y\mid\Theta=\theta)\sim1/\log(\norm{\theta}^2/\tau^2)$.
The weighted density $\pi_{p,\tau}(\theta)E(\kappa_Y\mid\Theta=\theta)$
is regularly varying, and the Gaussian convolution argument of
Appendix~\ref{sec:large_signal}, which extends unchanged to spherically
symmetric densities, gives
\begin{equation}
E(\kappa_Y\mid D=d)
\sim\frac{1}{\log(\norm{d}^2/\tau^2)},
\qquad \norm{d}\to\infty.
\label{eq:caucau_latent_asymptotic}
\end{equation}
For the calibrated coordinate, with $\varphi=\tau^2/\sigma^2$ as in
Section~\ref{sec:kappa_posterior}, the same argument gives
\begin{equation}
E(\kappa_Y^\star\mid D=d)
\sim\frac{A(\varphi)}{\log(\norm{d}^2/\tau^2)},
\qquad
A(\varphi)=\frac{\log\varphi}{\varphi-1},\quad A(1)=1,
\label{eq:caucau_latent_calibrated_asymptotic}
\end{equation}
since
$\int_0^\infty v\{(1+v^2)(\varphi+v^2)\}^{-1}dv=A(\varphi)/2$.
By contrast, \eqref{eq:caucau_mv_mean_tail} gives
\begin{equation}
\kappa_{\rm eff}(d)=1-a_p(\norm{d})
\sim\frac{(p+1)\sigma^2}{\norm{d}^2},
\qquad \norm{d}\to\infty.
\label{eq:caucau_eff_asymptotic}
\end{equation}

Latent denominator contraction and effective posterior shrinkage
therefore differ in rate, not merely in constant: the denominator
coordinate decays logarithmically, whereas posterior-mean attenuation
decays quadratically. The Cauchy denominator retains the exact horseshoe Beta
geometry, but the heavy-tailed numerator changes its relationship to
posterior shrinkage. This contrasts with the finite-moment setting of
Theorem~\ref{thm:kappa_large_signal}, where the two quantities share
the $d^{-2}$ rate.

\section{Illustration: Wavelet Denoising of Standard Test Signals}
\label{sec:dj_illustration}

The preceding analysis establishes the structural and asymptotic
properties of the Cauchy/Cauchy prior. To examine the finite-sample
behavior of its scalar posterior mean as a coefficientwise wavelet
shrinkage rule, we apply it to the four standard test signals of
\citet{DonohoJohnstone1994}: Blocks, Bumps, Doppler, and HeaviSine.
Each signal is sampled at $t_i=i/(N+1)$, $i=1,\ldots,N$, with $N=1024$,
and rescaled so that the root signal-to-noise ratio
${\rm sd}(f)/\sigma$ equals $5$, where ${\rm sd}(f)$ is the
population standard deviation of the sampled signal. Observations are
$y_i=f(t_i)+\varepsilon_i$ with independent $N(0,1)$ errors. The noise
level is used only to generate the data. In every replication it is
estimated from the finest-level detail coefficients $d_{J-1,\cdot}$ by
$\widehat\sigma={\rm MAD}(d_{J-1,\cdot})/\Phi^{-1}(3/4)$.
This estimate is supplied to methods that require an externally
estimated Gaussian noise scale. We use periodic orthonormal wavelet transforms matched to
the signal morphology: Haar for Blocks, Daubechies' six-tap filter for
Bumps, and the eight-tap Symmlet for Doppler and HeaviSine
\citep{daubechies1992ten}. The eight scaling coefficients at level
$j_0=3$ are kept, and detail levels $j=3,\ldots,9$ are shrunk. All
methods are applied to the same $M=1000$ noisy realizations of each
signal, so comparisons between methods are paired.

\paragraph{Cauchy/Cauchy rule.}
At level $j$ the coefficients are treated as independent observations
$d_{jk}\mid\theta_{jk}\sim N(\theta_{jk},\sigma^2)$ with prior
\eqref{eq:caucau_density} and scale $\tau_j$. Writing
$|\theta|=\tau_je^{t}$, the variable $t=\log(|\theta|/\tau_j)$ has density
$\ell(t)=2t/\{\pi^2\sinh t\}$ on $\R$, so for $z=|d|/\widehat\sigma$ and
$a(t)=\tau_je^{t}/\widehat\sigma$,
$$
\widehat\theta_{\rm mean}(d)=
{\rm sign}(d)\,\widehat\sigma\,
\frac{\int a(t)\,\ell(t)\{\phi(z-a(t))-\phi(z+a(t))\}\,dt}
     {\int \ell(t)\{\phi(z-a(t))+\phi(z+a(t))\}\,dt},
$$
a one-dimensional integral that we evaluate by the trapezoidal rule on
a fine grid in $t$. We report two choices of $\tau_j$. The first is a
fixed geometric schedule $\tau_j=45\,\widehat\sigma\,2^{-1.45j}$, chosen
once in exploratory runs under a different signal-to-noise convention
and held fixed for all four signals. The second contains no constants tuned by mean squared error. Following
the marginal maximum likelihood approach of
\citet{vanderPasSzaboVaart2017}, adapted to one scale per resolution
level, we set
$$
\widehat\tau_j=\arg\max_{\tau}\prod_k m_{\tau}(d_{jk}),
$$
where $m_\tau$ is the marginal density of $d_{jk}$ under the
Cauchy/Cauchy prior with scale $\tau$ and noise level $\widehat\sigma$.
The maximum is taken over a logarithmic grid on
$[\widehat\sigma/n_j,\,16\widehat\sigma]$, where $n_j=2^j$ is the number
of coefficients at level $j$.

\paragraph{Competing rules.}
We compare with fifteen established procedures, most of them Bayesian
or empirical Bayes:
the horseshoe posterior mean with $\tau_j$ chosen by marginal maximum
likelihood on $[\widehat\sigma/n_j,\widehat\sigma]$
\citep{carvalho2010horseshoe,vanderPasSzaboVaart2017};
BAMS \citep{VidakovicRuggeri2001};
the empirical Bayes posterior median and mean under a Laplace slab with
scale parameter $a=0.5$ and mixing weight by marginal maximum likelihood
\citep{JohnstoneSilverman2005};
BayesThresh, the posterior median under a normal-mixture prior with
hyperparameters $(0.5,1)$ in the $(\alpha,\beta)$ notation of that paper and constants $(C_1,C_2)$ by marginal maximum
likelihood \citep{AbramovichSapatinasSilverman1998};
raised-cosine \citep{ReinaSousa2026} and Epanechnikov
\citep{CausilSousa2025} Bayes rules;
the noncentral-$\chi^2$ block posterior mean NCPmean with blocks of
size two \citep{WangWood2006};
$\Lambda$-neighborhood shrinkage, LNWS \citep{RemenyiVidakovic2013};
BayesShrink \citep{ChangYuVetterli2000};
the bivariate parent--child rule BiShrink \citep{SendurSelesnick2002};
and four non-Bayesian benchmarks: SureShrink
\citep{DonohoJohnstone1995}, smooth SCAD with SURE-selected threshold
\citep{KulkarniEtAl2026}, BlockJS \citep{Cai1999} and VisuShrink
\citep{DonohoJohnstone1994}. BAMS, NCPmean, LNWS and the Epanechnikov
rule retain their own internal noise-scale estimates. No method uses
the true value of $\sigma$.

\begin{table}[thb]
\centering
\caption{Average mean squared error (AMSE) over $M=1000$ replications,
$N=1024$, root SNR $=5$, estimated $\sigma$. Ranks among the seventeen
rules are in parentheses; the smallest AMSE for each signal is in bold.
Monte Carlo standard errors of the AMSEs are at most $0.0023$ (at most
$0.0048$ for VisuShrink). Rules other than the two Cauchy/Cauchy
variants are ordered by mean rank. $^{\ast}$Non-Bayesian rule.}
\label{tab:dj_amse}
\small
\setlength{\tabcolsep}{3.2pt}
\begin{tabular}{@{}l r@{\,}l r@{\,}l r@{\,}l r@{\,}l c@{}}
\toprule
Rule & \multicolumn{2}{c}{Blocks} & \multicolumn{2}{c}{Bumps}
     & \multicolumn{2}{c}{Doppler} & \multicolumn{2}{c}{HeaviSine}
     & Mean rank\\
\midrule

Cau/Cau posterior mean, MMLE $\tau_j$ & 0.1196 & (5) & 0.2923 & (2) & 0.1652 & (6) & \textbf{0.0549} & (1) & 3.50 \\
Cau/Cau posterior mean, fixed $\tau_j$ & 0.1109 & (2) & 0.3209 & (9) & 0.1620 & (4) & 0.0623 & (10) & 6.25 \\
\midrule
BAMS & \textbf{0.1099} & (1) & 0.3270 & (10) & 0.1643 & (5) & 0.0575 & (4) & 5.00 \\
Horseshoe posterior mean, MMLE $\tau_j$ & 0.1248 & (7) & 0.2997 & (5) & 0.1689 & (7) & 0.0581 & (5) & 6.00 \\
NCPmean ($m=2$) & 0.1610 & (12) & \textbf{0.2727} & (1) & 0.1473 & (3) & 0.0606 & (8) & 6.00 \\
$\Lambda$-neighborhood (LNWS) & 0.1330 & (10) & 0.2931 & (4) & \textbf{0.1250} & (1) & 0.0701 & (12) & 6.75 \\
EBayesThresh, posterior median & 0.1233 & (6) & 0.3139 & (7) & 0.1761 & (10) & 0.0591 & (6) & 7.25 \\
BayesThresh, posterior median & 0.1147 & (3) & 0.3353 & (11) & 0.1946 & (13) & 0.0558 & (2) & 7.25 \\
EBayesThresh, posterior mean & 0.1329 & (9) & 0.3062 & (6) & 0.1702 & (8) & 0.0595 & (7) & 7.50 \\
Bivariate shrinkage (BiShrink) & 0.1982 & (14) & 0.2926 & (3) & 0.1257 & (2) & 0.0846 & (15) & 8.50 \\
Raised-cosine Bayes & 0.1166 & (4) & 0.3800 & (12) & 0.2021 & (15) & 0.0615 & (9) & 10.00 \\
Smooth SCAD--SURE$^{\ast}$ & 0.1257 & (8) & 0.3141 & (8) & 0.1852 & (11) & 0.0740 & (13) & 10.00 \\
Epanechnikov Bayes & 0.1497 & (11) & 0.5089 & (15) & 0.1862 & (12) & 0.0572 & (3) & 10.25 \\
SureShrink$^{\ast}$ & 0.1956 & (13) & 0.3983 & (13) & 0.2011 & (14) & 0.0821 & (14) & 13.50 \\
BayesShrink & 0.2948 & (15) & 0.4335 & (14) & 0.2331 & (16) & 0.0633 & (11) & 14.00 \\
BlockJS$^{\ast}$ & 0.4330 & (16) & 0.5977 & (16) & 0.1752 & (9) & 0.1137 & (16) & 14.25 \\
VisuShrink$^{\ast}$ & 0.7252 & (17) & 1.7262 & (17) & 0.6468 & (17) & 0.1985 & (17) & 17.00 \\
\bottomrule
\end{tabular}
\end{table}

\begin{figure}[t]
\centering
\includegraphics[width=\textwidth]{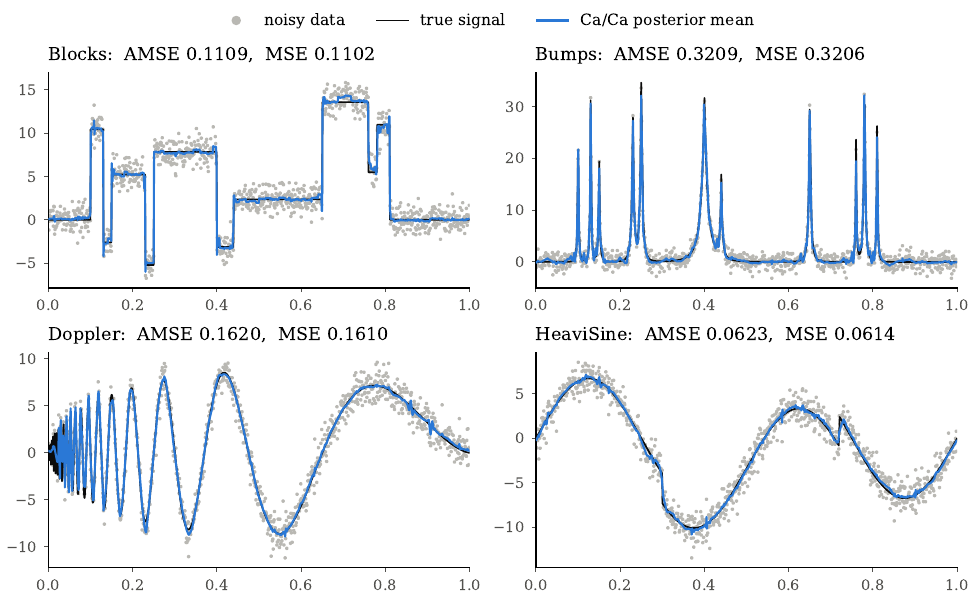}
\caption{One replication of each test signal (root SNR $=5$,
$N=1024$): noisy observations, the true signal, and the Cauchy/Cauchy
posterior mean with the fixed schedule
$\tau_j=45\,\widehat\sigma\,2^{-1.45j}$. The replication shown is the one
with the median MSE. Panel titles give the AMSE over $M=1000$
replications (as in Table~\ref{tab:dj_amse}) and the MSE of the
replication shown.}
\label{fig:dj_examples}
\end{figure}

\paragraph{Results.}
Table~\ref{tab:dj_amse} reports the AMSEs and
Figure~\ref{fig:dj_examples} shows typical reconstructions. No single
rule is best for all four signals; the per-signal winners are BAMS,
NCPmean, LNWS and the Cauchy/Cauchy rule with marginal-likelihood
hyperparameters. With this data-driven choice of $\tau_j$, the Cauchy/Cauchy
posterior mean has the best mean rank among the seventeen rules. It
ranks first for HeaviSine, second for Bumps, and fifth and sixth for
Blocks and Doppler. The fixed schedule is better for Blocks, where it
ranks second, but worse for Bumps and HeaviSine.

The most direct comparison is with the horseshoe, because both rules
are posterior means under heavy-tailed priors with
marginal-likelihood scales, applied coefficientwise. The Cauchy/Cauchy
rule has the smaller AMSE for every signal, with paired differences of
$-0.0052$, $-0.0074$, $-0.0036$ and $-0.0032$ (each with standard error
at most $0.0001$), and the smaller replication MSE in $99\%$, $99\%$,
$92\%$ and $91\%$ of the replications. Widening the horseshoe's range
of $\tau_j$ to $[\widehat\sigma/n_j,16\widehat\sigma]$, the range used
for the Cauchy/Cauchy rule, did not change this ordering.

On Doppler, the three rules that pool information across coefficients
(LNWS, BiShrink and NCPmean) have the three smallest AMSEs, and NCPmean
also has the smallest AMSE on Bumps. Both Cauchy/Cauchy variants
shrink each coefficient separately; the multivariate prior of
Section~\ref{sec:caucau_prior}, which shrinks a block through its
norm, is a natural way to add such pooling.

Finally, the estimated scales
$\widehat\tau_j$ decrease roughly geometrically with $j$, and for
Blocks their median rate is close to that of the fixed schedule. This
suggests that a prescribed geometric decay is a reasonable default
when marginal-likelihood estimation is not used.

\section{Discussion and Conclusions}
\label{sec:discussion}

The ratio representation provides a unified construction of shrinkage
priors in which the numerator and denominator play distinct but
complementary roles. The numerator determines the directional geometry
and contributes to scale and tail behavior; large values of the
denominator generate concentration near zero, and small values allow
escape into the tails. Through regular variation, our results
characterize both mechanisms, including power and logarithmic origin
singularities, their dependence on dimension, and the influence of
heavy-tailed numerators on the induced tails.

These tails have a direct inferential consequence. Under a Gaussian
likelihood, smoothly regularly varying prior tails yield posterior
means that approach the observation for large signals, with
squared-error risk approaching $p\sigma^2$. Tail properties of
ratio-generated priors thus translate directly into large-signal
robustness.

The denominator also induces a shrinkage geometry through
$\kappa_Y=Y^2/(1+Y^2)$. For a standard Cauchy denominator,
$\kappa_Y\sim{\rm Beta}(1/2,1/2)$.
With a Gaussian numerator, the horseshoe shrinkage coefficient is
the calibrated coordinate $\kappa_Y^\star$ in
\eqref{eq:kappaY_star}, which equals $\kappa_Y$ when $\tau=\sigma$.
With non-Gaussian numerators, latent denominator contraction and
effective posterior shrinkage can differ substantially: in the
Cauchy/Cauchy example they decay at different asymptotic rates.
Moreover, distinct ratio representations can generate the same
marginal prior with different latent geometries. Horseshoe-type
marginal behavior therefore does not determine the latent shrinkage
geometry, and that geometry does not necessarily determine effective
posterior shrinkage.

Because the ratio construction specifies a prior rather than an
estimator, different posterior summaries can produce different
shrinkage rules. For the Cauchy/Cauchy prior, the posterior mean
shrinks smoothly, whereas the posterior spatial median behaves nearly
as a threshold rule in our numerical illustrations. Exact variable
selection requires a decision rule, such as thresholding or posterior classification based on a practical null
region. For the Cauchy/Cauchy prior, one can also base a selection
rule on $E(\kappa_{VY}\mid D=d)$ in \eqref{eq:caucau_kappa_star}.

The present theory concerns proper ratio-generated priors with
positive regular-variation indices. The normal/log-Cauchy construction
of Remark~\ref{rem:logcauchy} shows what happens at the boundary
$\alpha=\beta=0$: the prior remains proper and close to the
scale-invariant density $1/|\theta|$, with a smaller large-signal
correction but only logarithmic control by the global scale. Extensions
based on improper numerator or denominator kernels are also possible,
although posterior propriety must then be established under the
specified likelihood.

Finally, the focus of this paper is theoretical. The wavelet
illustration of Section~\ref{sec:dj_illustration} shows that the
univariate Cauchy/Cauchy posterior mean, with level scales chosen by
marginal maximum likelihood, is competitive with established Bayesian
and non-Bayesian wavelet shrinkage rules. It has the best average rank
among seventeen rules on the four Donoho--Johnstone test signals, and
it improves on the horseshoe posterior mean computed the same way for
every signal. A systematic study is left for future work. Such a study
should cover other signal-to-noise ratios and sample sizes, the
multivariate (block) Cauchy/Cauchy rule, posterior contraction under
ML II, and variable selection in linear
regression.

\paragraph{Code availability.}
The code and simulation outputs for Table~\ref{tab:dj_amse} and Figure~\ref{fig:dj_examples} are
available at \url{https://github.com/BraniV/HorseshoeTypeRatio}.
The repository includes results from all 1,000 replications for
each of the four test signals and instructions for rerunning the
simulations and reproducing the table and figure.

\section*{Appendix}
\addcontentsline{toc}{section}{Appendix}

\renewcommand{\thesubsection}{A.\arabic{subsection}}
\setcounter{subsection}{0}

 
\subsection{Proof of Theorem~\ref{thm:origin_power}}
\label{sec:prooftheorem2}
 
\begin{proof}
Fix a unit vector $e$ and $\eta\in(0,p-\alpha)$ as in the assumptions,
and let $B=\sup_{0\leq u\leq1}f(ue)<\infty$.
 
By \eqref{eq:q_origin_rv}, $q$ is regularly varying with index
$-1-\alpha$ and $q(1/r)\sim c_Or^{1+\alpha}L_O(1/r)$ as $r\downarrow0$.
Potter's bounds \citep[Theorem~1.5.6]{BinghamGoldieTeugels1987},
applied to $q(u/r)/q(1/r)$, therefore give $y_0\geq1$, $r_0>0$ and
$K_\eta<\infty$ such that, for $0<r<r_0$ and $u\geq ry_0$,
$$
\frac{q(u/r)}{c_Or^{1+\alpha}L_O(1/r)}
\leq K_\eta\left(u^{-1-\alpha-\eta}+u^{-1-\alpha+\eta}\right).
$$
Moreover, for every fixed $u>0$,
$$
\frac{q(u/r)}{c_Or^{1+\alpha}L_O(1/r)}\longrightarrow u^{-1-\alpha},
\qquad r\downarrow0.
$$
Reducing $r_0$ if necessary, assume that $r_0y_0\leq1$.
 
Split \eqref{eq:radial_ratio_density} as $\pi(re)=I_1(r)+I_2(r)$, where
$$
I_1(r)=\int_0^{y_0}y^pf(rye)q(y)\,dy,
\qquad
I_2(r)=\int_{y_0}^\infty y^pf(rye)q(y)\,dy.
$$
For $y\leq y_0$ we have $ry\leq1$, so $f(rye)\leq B$, and since $q$ is a
density,
$$
0\leq I_1(r)\leq B\int_0^{y_0}y^pq(y)\,dy\leq By_0^p.
$$
Because $\alpha<p$, $r^{\alpha-p}L_O(1/r)\to\infty$ as $r\downarrow0$,
and therefore $I_1(r)=o\{r^{\alpha-p}L_O(1/r)\}$.
 
For $I_2(r)$, the substitution $u=ry$ gives
$$
\frac{I_2(r)}{c_Or^{\alpha-p}L_O(1/r)}
=\int_0^\infty u^pf(ue)Q_r(u)\,du,
\qquad
Q_r(u)=\mathbf{1}_{\{u\geq ry_0\}}
\frac{q(u/r)}{c_Or^{1+\alpha}L_O(1/r)}.
$$
By the Potter bound,
$$
0\leq u^pf(ue)Q_r(u)
\leq K_\eta\left(u^{p-\alpha-1-\eta}+u^{p-\alpha-1+\eta}\right)f(ue).
$$
The right-hand side is integrable on $(0,1]$, because $f(ue)\leq B$
there and $p-\alpha-1-\eta>-1$, and on $[1,\infty)$, because
$u^{p-\alpha-1-\eta}\leq u^{p-\alpha-1+\eta}$ there and
$\int_1^\infty u^{p-\alpha-1+\eta}f(ue)\,du<\infty$. Since
$Q_r(u)\to u^{-1-\alpha}$ for every $u>0$, dominated convergence yields
$$
\lim_{r\downarrow0}\frac{I_2(r)}{c_Or^{\alpha-p}L_O(1/r)}
=\int_0^\infty u^{p-\alpha-1}f(ue)\,du=C_\alpha(e).
$$
Combining this with $I_1(r)=o\{r^{\alpha-p}L_O(1/r)\}$ establishes
\eqref{eq:origin_power}.
\end{proof}
 
 
\subsection{Proof of Proposition~\ref{prop:origin_log}}
\label{sec:proofproposition}
 
\begin{proof}
Fix $\varepsilon\in(0,1)$. By continuity of $f$ at the origin, choose
$\delta\in(0,1)$ such that
$$
(1-\varepsilon)f(0)\leq f(ue)\leq(1+\varepsilon)f(0),
\qquad 0\leq u\leq\delta,
$$
and by \eqref{eq:q_origin_rv} with $\alpha=p$, choose $y_0$ large
enough that
$$
(1-\varepsilon)c_O\frac{L_O(y)}{y^{p+1}}
\leq q(y)\leq
(1+\varepsilon)c_O\frac{L_O(y)}{y^{p+1}},
\qquad y\geq y_0.
$$
For $r<\delta/y_0$, split \eqref{eq:radial_ratio_density} as
$\pi(re)=I_1(r)+I_2(r)+I_3(r)$, with integration over $(0,y_0)$,
$(y_0,\delta/r)$ and $(\delta/r,\infty)$, respectively.
 
For $y<y_0$ we have $ry<\delta$, so $f(rye)\leq(1+\varepsilon)f(0)$ and
$$
0\leq I_1(r)\leq(1+\varepsilon)f(0)\int_0^{y_0}y^pq(y)\,dy
\leq(1+\varepsilon)f(0)y_0^p.
$$
Hence $I_1(r)=O(1)=o\{\widetilde L_O(1/r)\}$.
 
For $y_0<y<\delta/r$, both two-sided bounds apply, so
$$
(1-\varepsilon)^2c_Of(0)\int_{y_0}^{\delta/r}\frac{L_O(y)}{y}\,dy
\leq I_2(r)\leq
(1+\varepsilon)^2c_Of(0)\int_{y_0}^{\delta/r}\frac{L_O(y)}{y}\,dy.
$$
By \citet[Proposition~1.5.9a]{BinghamGoldieTeugels1987}, $\widetilde L_O$
is slowly varying and $L_O(x)=o\{\widetilde L_O(x)\}$. Since
$\widetilde L_O(x)\to\infty$,
$$
\int_{y_0}^{\delta/r}\frac{L_O(y)}{y}\,dy
=\widetilde L_O(\delta/r)-\widetilde L_O(y_0)
\sim\widetilde L_O(1/r).
$$
 
For $I_3(r)$, the substitution $u=ry$ gives
$$
I_3(r)=r^{-(p+1)}\int_\delta^\infty u^pf(ue)q(u/r)\,du.
$$
Since $q$ is regularly varying with index $-(p+1)$, Potter's bounds
give $A_\eta<\infty$ such that, for all sufficiently small $r$,
$$
\frac{q(u/r)}{q(1/r)}\leq A_\eta u^{-(p+1)}\max\{u^\eta,u^{-\eta}\},
\qquad u\geq\delta.
$$
Consequently,
$$
I_3(r)\leq A_\eta\frac{q(1/r)}{r^{p+1}}
\int_\delta^\infty\frac{f(ue)}{u}\max\{u^\eta,u^{-\eta}\}\,du.
$$
The integral is finite, by local boundedness of $f(ue)$ on
$[\delta,1]$ and by \eqref{eq:boundary_ray_integrability} on
$[1,\infty)$. Since $q(1/r)/r^{p+1}\sim c_OL_O(1/r)$, we obtain
$I_3(r)=O\{L_O(1/r)\}=o\{\widetilde L_O(1/r)\}$.
 
Thus only $I_2(r)$ contributes to the leading order. Dividing its bounds
by $c_Of(0)\widetilde L_O(1/r)$, letting $r\downarrow0$ and then
$\varepsilon\downarrow0$ proves \eqref{eq:origin_log_general}. Finally,
if $L_O(y)\to L_{O,0}\in(0,\infty)$, then
$\widetilde L_O(x)\sim L_{O,0}\log x$, which gives \eqref{eq:origin_log}.
\end{proof}

\subsection{Proof of Theorem~\ref{thm:radial_tail}}
\label{sec:proofth3}

\begin{proof}
The argument is the standard proof of Breiman's lemma
\citep{Breiman1965,ClineSamorodnitsky1994}, written in the present
notation for completeness.

Write $W=\norm{X}$ and choose $\eta\in(0,\min\{\beta,\epsilon\})$.
By \eqref{eq:H_regular_variation}, there is $t_0>0$ such that
$H(t)\leq2c_Yt^\beta L_T(1/t)$ for $0<t\leq t_0$. Potter's bounds
\citep[Theorem~1.5.6]{BinghamGoldieTeugels1987} give $A_\eta<\infty$
and $x_0$ such that $L_T(x)/L_T(y)\leq A_\eta\max\{(x/y)^\eta,(x/y)^{-\eta}\}$
for $x,y\geq x_0$. Reducing $t_0$ so that $1/t_0\geq x_0$, we obtain,
for all $r\geq x_0$ and $0<W\leq rt_0$,
$$
\frac{H(W/r)}{r^{-\beta}L_T(r)}
\leq 2c_YW^\beta\,\frac{L_T(r/W)}{L_T(r)}
\leq C_\eta\left(W^{\beta-\eta}+W^{\beta+\eta}\right),
$$
with $C_\eta=2c_YA_\eta$ independent of $r$ and $W$. The dominating
function is integrable, since $\beta-\eta>0$ and
$EW^{\beta+\epsilon}<\infty$ with $\eta<\epsilon$.

For each fixed $w>0$, regular variation gives
$$
\frac{H(w/r)}{r^{-\beta}L_T(r)}
=c_Yw^\beta\,\frac{L_T(r/w)}{L_T(r)}\{1+o(1)\}
\longrightarrow c_Yw^\beta,
$$
and the ratio is identically zero at $w=0$, since $H(0)=0$. Because
$\mathbf{1}_{\{W\leq rt_0\}}\to1$ almost surely, dominated convergence
gives
$$
\lim_{r\to\infty}
\frac{E\{H(W/r);\,W\leq rt_0\}}{r^{-\beta}L_T(r)}
=c_YEW^\beta.
$$

For the complementary contribution, $H\leq1$ and Markov's inequality
yield
$$
0\leq E\{H(W/r);\,W>rt_0\}
\leq P(W>rt_0)
\leq\frac{EW^{\beta+\epsilon}}{(rt_0)^{\beta+\epsilon}}
=o\{r^{-\beta}L_T(r)\},
$$
since $r^{-\epsilon}/L_T(r)\to0$ for slowly varying $L_T$. Combining
the two parts with the identity $P(\norm{\Theta}>r)=E\{H(W/r)\}$ of
\eqref{eq:tail_identity} establishes \eqref{eq:radial_tail}.
\end{proof}

\subsection{Proof of Theorem~\ref{thm:density_tail}}
\label{sec:proofth4}
\begin{proof}
Fix $\delta>0$ sufficiently small and split
\eqref{eq:radial_ratio_density} as
$\pi(re)=I_1(r)+I_2(r)$, with integration over
$(0,\delta)$ and $(\delta,\infty)$, respectively.

For $I_1(r)$, substituting $u=ry$ gives
$$
\frac{I_1(r)}
{c_Tr^{-(p+\beta)}L_T(r)}
=
\int_0^{r\delta}
u^pf(ue)
\frac{q(u/r)}
{c_Tr^{1-\beta}L_T(r)}\,du.
$$
By regular variation in \eqref{eq:q_tail_rv},
the normalized integrand, extended by zero for
$u>r\delta$, converges pointwise to
$u^{p+\beta-1}f(ue)$.
Potter's bounds dominate it by
$$
C_\eta
\left(
u^{p+\beta-1-\eta}
+
u^{p+\beta-1+\eta}
\right)f(ue),
$$
which is integrable by assumption. Dominated
convergence therefore yields
$$
\frac{I_1(r)}
{c_Tr^{-(p+\beta)}L_T(r)}
\longrightarrow C_\beta(e).
$$

For $I_2(r)$, boundedness of $q$ on
$[\delta,\infty)$ gives, for sufficiently large $r$,
$$
\begin{aligned}
I_2(r)
&\leq C_\delta r^{-(p+1)}
\int_{r\delta}^\infty u^pf(ue)\,du\\
&\leq C_\delta r^{-(p+1)}
(r\delta)^{-(\beta-1+\eta)}
\int_{r\delta}^\infty
u^{p+\beta-1+\eta}f(ue)\,du\\
&=O\left(r^{-(p+\beta+\eta)}\right).
\end{aligned}
$$
Here $\beta-1+\eta>0$, and the weighted integral
is finite by assumption. Since $L_T$ is positive
and slowly varying,
$$
\frac{I_2(r)}
{r^{-(p+\beta)}L_T(r)}
=
O\left(\frac{r^{-\eta}}{L_T(r)}\right)
\longrightarrow0.
$$
Thus $I_2(r)$ is negligible relative to $I_1(r)$,
and \eqref{eq:density_tail} follows.
\end{proof} 

\subsection{Proof of Lemma~\ref{lem:gaussian_convolution_tail}}
\label{sec:prooflemma}

\begin{proof}
Both $\pi$ and $\phi_{\sigma,p}$ are spherically symmetric, so $m$ is
spherically symmetric, and it suffices to consider $d=te_1$ with
$t\to\infty$, where $e_1=(1,0,\ldots,0)^\top\in\R^p$. Write
$\rho=p+\beta$. Since $\phi_{\sigma,p}$ is symmetric, the substitution
$u=\theta-d$ gives
$$
m(d)=\int_{\R^p}\phi_{\sigma,p}(u)\pi(d+u)\,du.
$$
Set
$$
A_t=\{u:\norm{u}\leq t/2\},
\qquad
B_t=\R^p\setminus A_t,
\qquad
R_t(u)=\frac{\pi(te_1+u)}{\pi(te_1)}.
$$

\emph{A uniform bound on $A_t$.} By
\eqref{eq:smooth_regular_variation}, there exist $t_0>0$ and
$C<\infty$ such that
\begin{equation}
\left|\frac{h'(r)}{h(r)}\right|\leq\frac{C}{r},
\qquad r\geq t_0.
\label{eq:proof_score_bound}
\end{equation}
Let $t\geq2t_0$ and $u\in A_t$. Every point $x=te_1+su$ with
$0\leq s\leq1$ satisfies $\norm{x}\geq t-\norm{u}\geq t/2\geq t_0$, and
$\norm{\nabla\log\pi(x)}=|h'(\norm{x})/h(\norm{x})|\leq2C/t$ there.
Integrating along the segment,
\begin{equation}
\log R_t(u)=\int_0^1\nabla\log\pi(te_1+su)^\top u\,ds,
\qquad
|\log R_t(u)|\leq\frac{2C\norm{u}}{t}\leq C.
\label{eq:proof_log_ratio}
\end{equation}
Consequently, uniformly over $u\in A_t$,
$$
R_t(u)\leq e^C,
\qquad
|R_t(u)-1|\leq e^C|\log R_t(u)|\leq C_3\frac{\norm{u}}{t},
\qquad C_3=2Ce^C.
$$

\emph{The marginal.} By \eqref{eq:proof_log_ratio}, $R_t(u)\to1$ for
each fixed $u$, and $0 \leq R_t\mathbf 1_{A_t}\leq e^C$. Dominated convergence
therefore yields
\begin{equation}
\int_{A_t}\phi_{\sigma,p}(u)R_t(u)\,du\longrightarrow1.
\label{eq:proof_central_mass}
\end{equation}
On $B_t$, since $\pi$ is a probability density,
$$
\int_{B_t}\phi_{\sigma,p}(u)\pi(te_1+u)\,du
\leq\sup_{\norm{u}>t/2}\phi_{\sigma,p}(u)
=(2\pi\sigma^2)^{-p/2}\exp\left\{-\frac{t^2}{8\sigma^2}\right\}.
$$
By \eqref{eq:prior_regular_tail_gaussian},
$\pi(te_1)=h(t)=t^{-\rho}L(t)\{1+o(1)\}$ decays more slowly than any
exponential, so this bound is $o\{\pi(te_1)\}$. Together with
\eqref{eq:proof_central_mass}, this establishes
\begin{equation}
m(d)\sim\pi(d),
\qquad \norm{d}\to\infty.
\label{eq:proof_marginal_equivalence}
\end{equation}

\emph{The score.} Condition \eqref{eq:smooth_regular_variation} gives
$$
\nabla\log\pi(x)
=\frac{h'(\norm{x})}{h(\norm{x})}\frac{x}{\norm{x}}
=-\rho\frac{x}{\norm{x}^2}+o(\norm{x}^{-1}),
\qquad \norm{x}\to\infty.
$$
For fixed $u$, $t\nabla\log\pi(te_1+su)\to-\rho e_1$ uniformly in
$s\in[0,1]$, so \eqref{eq:proof_log_ratio} gives
$t\log R_t(u)\to-\rho u_1$, where $u_1$ is the first coordinate of $u$.
Since $\log R_t(u)\to0$, it follows that
\begin{equation}
t\{R_t(u)-1\}\longrightarrow-\rho u_1.
\label{eq:proof_first_order}
\end{equation}
Because $A_t$ is symmetric about the origin,
$\int_{A_t}u\,\phi_{\sigma,p}(u)\,du=0$, and therefore
$$
\frac{t}{\pi(te_1)}\int_{A_t}u\,\phi_{\sigma,p}(u)\pi(te_1+u)\,du
=\int_{A_t}u\,\phi_{\sigma,p}(u)\,t\{R_t(u)-1\}\,du.
$$
The norm of the integrand is at most $C_3\norm{u}^2\phi_{\sigma,p}(u)$,
which is integrable. By \eqref{eq:proof_first_order} and dominated
convergence,
\begin{equation}
\frac{t}{\pi(te_1)}\int_{A_t}u\,\phi_{\sigma,p}(u)\pi(te_1+u)\,du
\longrightarrow
-\rho\int_{\R^p}u\,u_1\phi_{\sigma,p}(u)\,du
=-\rho\sigma^2e_1.
\label{eq:proof_score_central}
\end{equation}
On $B_t$, the function $r\mapsto re^{-r^2/(2\sigma^2)}$ is decreasing for
$r>\sigma$, so for $t\geq2\sigma$,
$$
\left\|\int_{B_t}u\,\phi_{\sigma,p}(u)\pi(te_1+u)\,du\right\|
\leq\sup_{\norm{u}>t/2}\norm{u}\phi_{\sigma,p}(u)
=\frac{t}{2}(2\pi\sigma^2)^{-p/2}\exp\left\{-\frac{t^2}{8\sigma^2}\right\},
$$
and hence
\begin{equation}
\frac{t}{\pi(te_1)}\int_{B_t}u\,\phi_{\sigma,p}(u)\pi(te_1+u)\,du
\longrightarrow0.
\label{eq:proof_score_remote}
\end{equation}
Combining \eqref{eq:proof_score_central} and
\eqref{eq:proof_score_remote},
$$
\int_{\R^p}u\,\phi_{\sigma,p}(u)\pi(te_1+u)\,du
=-\frac{\rho\sigma^2}{t}\pi(te_1)e_1+o\left\{\frac{\pi(te_1)}{t}\right\}.
$$
By \eqref{eq:marginal_derivative}, with the substitution $u=\theta-d$,
the left-hand side equals $\sigma^2\nabla m(te_1)$. Dividing by
$\sigma^2m(te_1)$ and using \eqref{eq:proof_marginal_equivalence},
$$
\nabla\log m(te_1)=-\frac{\rho}{t}e_1+o(t^{-1}),
$$
and spherical symmetry of $m$ gives \eqref{eq:marginal_score_tail}.
\end{proof}


\subsection{Proof of Theorem~\ref{thm:posterior_mean_robustness}}
\label{sec:proofth5}

\begin{proof}
The expansion \eqref{eq:posterior_mean_expansion} follows from
\eqref{eq:brown_identity} and \eqref{eq:marginal_score_tail}. For the
risk, write
$$
a(d)=\widehat\theta_{\rm mean}(d)-d=\sigma^2\nabla\log m(d),
$$
by \eqref{eq:brown_identity}. The Gaussian convolution $m$ is strictly
positive and infinitely differentiable, so $a$ is continuous on $\R^p$,
and by \eqref{eq:posterior_mean_expansion}, $\norm{a(d)}\to0$ as
$\norm{d}\to\infty$. A continuous function that vanishes at infinity is
bounded, so
$$
M=\sup_{d\in\R^p}\norm{a(d)}<\infty.
$$

Write $D=\theta+\sigma Z$ with $Z\sim N_p(0,I_p)$. Then
$$
\widehat\theta_{\rm mean}(D)-\theta=\sigma Z+a(D).
$$
Let $\norm{\theta_n}\to\infty$. For every fixed $z$,
$\norm{\theta_n+\sigma z}\to\infty$, so $a(\theta_n+\sigma z)\to0$.
Since $\norm{a}^2\leq M^2$, dominated convergence gives
$$
E_{\theta_n}\norm{a(D)}^2\longrightarrow0.
$$
By the triangle inequality in $L^2$,
$$
\left|
R(\theta_n,\widehat\theta_{\rm mean})^{1/2}
-\left(E\norm{\sigma Z}^2\right)^{1/2}
\right|
\leq\left(E_{\theta_n}\norm{a(D)}^2\right)^{1/2}
\longrightarrow0,
$$
and $E\norm{\sigma Z}^2=p\sigma^2$. Since the sequence $(\theta_n)$ was
arbitrary, $R(\theta,\widehat\theta_{\rm mean})\to p\sigma^2$ as
$\norm{\theta}\to\infty$.
\end{proof}


\subsection{Proof of Theorem~\ref{thm:kappa_large_signal}}
\label{sec:large_signal}

\begin{proof}
By symmetry it suffices to take $\theta>0$. Let $U=|Z|$, with density
$f_U(u)=2f_Z(u)$ on $(0,\infty)$, put $t=\theta/\tau$, and write $P_t$
and $E_t$ for probability and expectation conditional on
$\Theta=\theta$.

\emph{Step 1: the conditional law of $U$.} Given $\Theta=\theta$, we
have $U=t|Y|$. The joint density of $(\Theta,Y)$ is
$|y|f_X(y\theta)g(y)$ and $f_X(x)=\tau^{-1}f_Z(x/\tau)$ is symmetric, so
the substitution $u=t|y|$ shows that, under $P_t$, $U$ has density
$$
p_t(u)=\frac{u f_U(u)q(u/t)}{N(t)},
\qquad
N(t)=\int_0^\infty v f_U(v)q(v/t)\,dv,
\qquad u>0.
$$
The same substitution in \eqref{eq:ratio_density} gives
\begin{equation}
\pi(\theta)=\frac{\tau}{2\theta^2}N(t).
\label{eq:proof_pi_N}
\end{equation}

\emph{Step 2: the normalizing constant.} Define
$$
n_t(u)=\frac{u f_U(u)q(u/t)}{c_Tt^{1-\beta}L_T(t)},
\qquad u>0,
$$
so that $p_t=n_t/\int_0^\infty n_t$. By \eqref{eq:q_tail_rv},
$n_t(u)\to u^\beta f_U(u)$ for every $u>0$. Fix
$\varepsilon\in(0,\min\{\beta,\eta\})$. By \eqref{eq:q_tail_rv} and
Potter's bounds for $L_T$, there are $\delta>0$ and $C<\infty$ such
that, for all sufficiently large $t$,
\begin{equation}
n_t(u)\leq C\left(u^{\beta-\varepsilon}+u^{\beta+\varepsilon}\right)
f_U(u),
\qquad 0<u\leq\delta t.
\label{eq:proof_potter_nt}
\end{equation}
Since $\varepsilon<\beta$, both powers are bounded near zero, and since
$\varepsilon<\eta$, the bound remains integrable after multiplication by
$u^2$, by the assumption $E|Z|^{\beta+2+\eta}<\infty$. On
$\{u>\delta t\}$, boundedness of $q$ on $[\delta,\infty)$ gives
\begin{equation}
\int_{\delta t}^\infty n_t(u)\,du
\leq\frac{\sup_{y\geq\delta}q(y)}{c_Tt^{1-\beta}L_T(t)}
\int_{\delta t}^\infty uf_U(u)\,du
=O\left(\frac{t^{-2-\eta}}{L_T(t)}\right),
\label{eq:proof_far_part}
\end{equation}
because
$\int_{\delta t}^\infty uf_U(u)\,du
\leq(\delta t)^{-\beta-1-\eta}E|Z|^{\beta+2+\eta}$. By dominated
convergence on $\{u\leq\delta t\}$ and \eqref{eq:proof_far_part},
$\int_0^\infty n_t\to\mu_\beta$, that is,
$$
N(t)\sim c_T\mu_\beta t^{1-\beta}L_T(t),
\qquad t\to\infty.
$$
By \eqref{eq:proof_pi_N}, with $t=\theta/\tau$,
\begin{equation}
\pi(\theta)\sim\frac{c_T\tau^\beta\mu_\beta}{2}\,
|\theta|^{-(1+\beta)}L_T(|\theta|/\tau),
\qquad|\theta|\to\infty,
\label{eq:proof_pi_tail}
\end{equation}
in agreement with Theorem~\ref{thm:density_tail} for $p=1$. In
particular, $\pi$ is regularly varying with index $-(1+\beta)$.

\emph{Step 3: the distributional limit.} Since
$p_t(u)=n_t(u)/\int_0^\infty n_t\to u^\beta f_U(u)/\mu_\beta$ for every
$u>0$, Scheff\'e's lemma gives convergence of $P_t(U\in\cdot)$ in total
variation 
(see \citealp[Section~2.2]{DasGupta2008})
to the law of $|Z_\beta|$. Since $Y=\tau Z/\theta$,
\begin{equation}
\theta^2\kappa_Y^\star
=\frac{\sigma^2\theta^2Z^2}{\theta^2+\sigma^2Z^2}
=\frac{\sigma^2U^2}{1+\sigma^2U^2/\theta^2}.
\label{eq:proof_kappa_U}
\end{equation}
Under $P_t$, $U$ is tight, so $\sigma^2U^2/\theta^2\to0$ in
probability, and Slutsky's lemma gives
$\theta^2\kappa_Y^\star\mid\{\Theta=\theta\}\Longrightarrow
\sigma^2Z_\beta^2$.

\emph{Step 4: the conditional expectation.} Split
$\theta^2E_t(\kappa_Y^\star)$ according to $\{U\leq\delta t\}$ and
$\{U>\delta t\}$. On the first set, \eqref{eq:proof_kappa_U} is bounded
by $\sigma^2U^2$, and \eqref{eq:proof_potter_nt} gives the integrable
dominating function
$C\sigma^2(u^{\beta+2-\varepsilon}+u^{\beta+2+\varepsilon})f_U(u)$; the
integrand converges pointwise to $\sigma^2u^{\beta+2}f_U(u)$. On the
second set, $\kappa_Y^\star\leq1$ and, by \eqref{eq:proof_far_part},
$$
\theta^2P_t(U>\delta t)
=O\left(\frac{t^2\cdot t^{-2-\eta}}{L_T(t)}\right)
=O\left(\frac{t^{-\eta}}{L_T(t)}\right)=o(1).
$$
Hence
\begin{equation}
\theta^2E(\kappa_Y^\star\mid\Theta=\theta)
\longrightarrow\sigma^2\frac{\mu_{\beta+2}}{\mu_\beta},
\qquad|\theta|\to\infty.
\label{eq:proof_cond_exp}
\end{equation}

\emph{Step 5: the posterior expectation.} Define
$w(\theta)=\pi(\theta)E(\kappa_Y^\star\mid\Theta=\theta)$. Since
$0\leq\kappa_Y^\star\leq1$, we have $0\leq w\leq\pi$, so $w$ is
integrable, and by \eqref{eq:proof_cond_exp},
$$
w(\theta)\sim\sigma^2\frac{\mu_{\beta+2}}{\mu_\beta}
\frac{\pi(\theta)}{\theta^2},
\qquad|\theta|\to\infty.
$$
By \eqref{eq:proof_pi_tail}, $w$ is regularly varying with index
$-(3+\beta)$. Writing $\phi_\sigma$ for the $N(0,\sigma^2)$ density,
$$
E(\kappa_Y^\star\mid D=d)
=\frac{(\phi_\sigma*w)(d)}{(\phi_\sigma*\pi)(d)}.
$$

For both convolutions we use the first part of the proof of
Lemma~\ref{lem:gaussian_convolution_tail}, which requires only regular
variation. Let $v\geq0$ be integrable on $\R$ and regularly varying with
negative index as $|x|\to\infty$, and let $d\to\infty$; the case
$d\to-\infty$ is symmetric. For $|u|\leq d/2$, we have
$(d+u)/d\in[1/2,3/2]$, so the uniform convergence theorem
\citep[Theorem~1.5.2]{BinghamGoldieTeugels1987} shows that
$v(d+u)/v(d)$ is bounded uniformly in $|u|\leq d/2$ and tends to one
for each fixed $u$. Dominated convergence gives
$\int_{|u|\leq d/2}\phi_\sigma(u)v(d+u)\,du\sim v(d)$, while
$$
\int_{|u|>d/2}\phi_\sigma(u)v(d+u)\,du
\leq\sup_{|u|>d/2}\phi_\sigma(u)\int_\R v
=O\left(e^{-d^2/(8\sigma^2)}\right)=o\{v(d)\}.
$$
Hence $(\phi_\sigma*v)(d)\sim v(d)$. Applying this with $v=w$ and
$v=\pi$,
$$
E(\kappa_Y^\star\mid D=d)\sim\frac{w(d)}{\pi(d)}
\sim\frac{\sigma^2}{d^2}\frac{\mu_{\beta+2}}{\mu_\beta},
$$
which proves \eqref{eq:kappa_latent_large}.

\emph{Step 6: the effective shrinkage.} By \eqref{eq:proof_pi_tail},
$\pi$ satisfies \eqref{eq:prior_regular_tail_gaussian} with $p=1$. If,
in addition, \eqref{eq:smooth_regular_variation} holds, then
Theorem~\ref{thm:posterior_mean_robustness} gives
$\widehat\theta_{\rm mean}(d)=d-\sigma^2(1+\beta)/d+o(d^{-1})$, which
is \eqref{eq:kappa_eff_large}. Dividing \eqref{eq:kappa_latent_large}
by \eqref{eq:kappa_eff_large} gives \eqref{eq:kappa_ratio_large}.
\end{proof}



\bibliographystyle{agsm}
\bibliography{references}

\end{document}